 \documentclass[final,1p,times]{elsarticleRossEdit}

\usepackage{graphicx, verbatim}
\usepackage{amssymb}
\usepackage{amsmath}

\newcommand{\bx}{{\mbox{\boldmath $x$}}}

\newcommand{\bu}{{\mbox{\boldmath $u$}}}

\newcommand{\bh}{{\mathbf h}}

\newcommand{\be}{{\mathbf e}}

\newcommand{\bnu}{{\mbox{\boldmath $\nu$}}}
\newcommand{\bmu}{{\mbox{\boldmath $\mu$}}}

\newcommand{\blam}{{\mbox{\boldmath $\lambda$}}}
\newcommand{\tlam}{{\mbox{$\widetilde{\lambda}$}}}

\newcommand{\tnu}{{\mbox{$\widetilde{\nu}$}}}

\newcommand{\abs}[1]{\left\vert#1\right\vert}
\newcommand{\set}[1]{\left\{#1\right\}}
\newcommand{\Real}{\mathbb R}
\newcommand{\real}[1]{{\mathbb R}^{#1}}
\newcommand{\U}{\mathbb U}
\newcommand{\X}{\mathbb X}

\newcommand{\bea}{\begin{eqnarray}}
\newcommand{\eea}{\end{eqnarray}}
\newcommand{\nn}{\nonumber}
\def\ds{\displaystyle}
\newcommand{\norm}[1]{\left\Vert#1\right\Vert}

\newtheorem{theorem}{Theorem}
\newtheorem{proposition}{Propostion}
\newtheorem{remark}{Remark}

\begin{document}

\begin{frontmatter}



\title{A Review of Pseudospectral Optimal Control: From Theory to Flight}


\author[mike]{I. Michael Ross}
\author[mark]{Mark Karpenko}

\address[mike]{Professor and Program Director, Control and Optimization, Naval Postgraduate School, Monterey, CA}
\address[mark]{Research Assistant Professor, Control and Optimization, Naval Postgraduate School, Monterey, CA}

\begin{abstract}
The home space for optimal control is a Sobolev space.  The home space for pseudospectral theory is also a Sobolev space.  It thus seems natural to combine pseudospectral theory with optimal control theory and construct ``pseudospectral optimal control theory,'' a term coined by Ross.
In this paper, we review key theoretical results in pseudospectral optimal control that have proven to be critical for a successful flight.  Implementation details of flight demonstrations onboard NASA spacecraft are discussed along with emerging trends and techniques in both theory and practice. The 2011 launch of pseudospectral optimal control in embedded platforms is changing the way in which we see solutions to challenging control problems in aerospace and autonomous systems.

\end{abstract}


\end{frontmatter}


\section{Introduction}
\label{intro}
\begin{figure}[b]
      \centering\scalebox{1}
      {\includegraphics[scale=0.35]{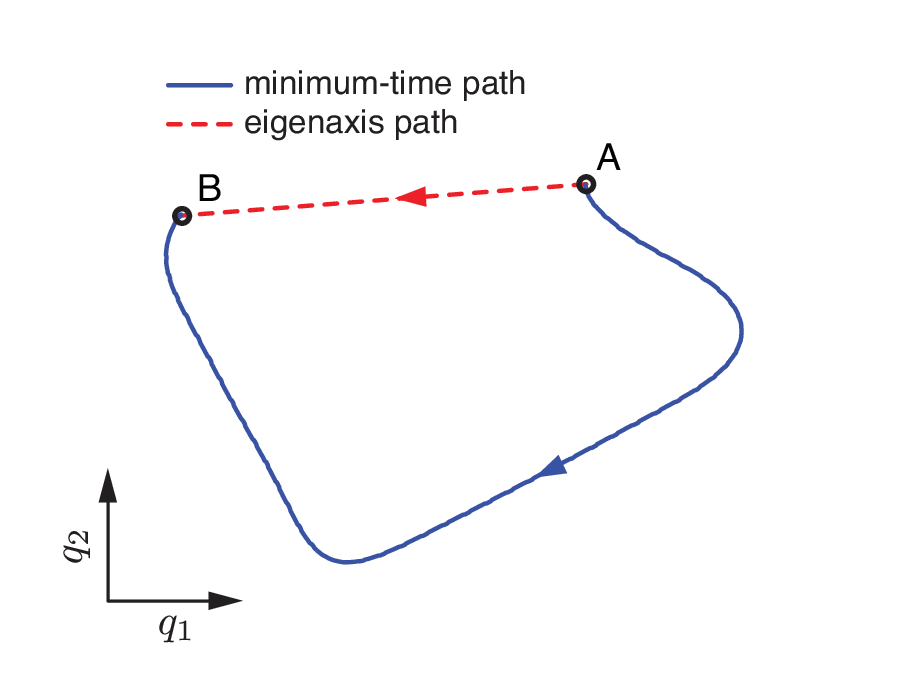}}
       {\includegraphics[scale=0.35]{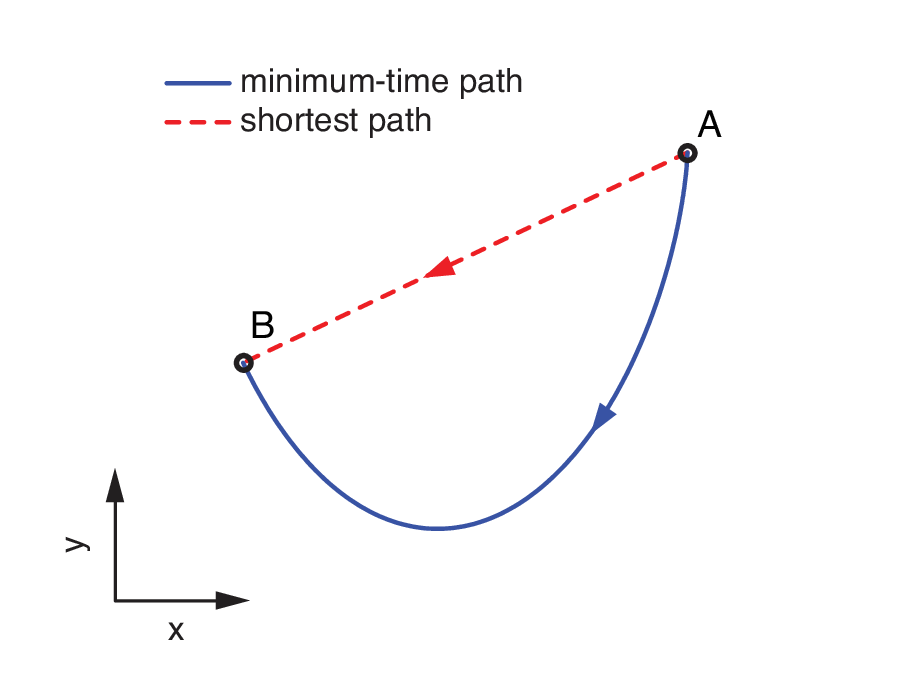}}
      \caption{\textsf{Flight results of TRACE's minimum-time maneuver (left) in projected ``quaternion space'' illustrating its analogy with the classical Brachistochrone problem (right).}}
      \label{fig:STM_Brac_q3q2}
   \end{figure}

\begin{figure}[h!]
      \centering
      {\includegraphics[scale=1.5]{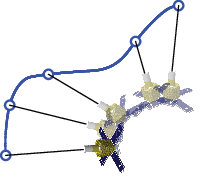}}
       {\includegraphics[scale=1.5]{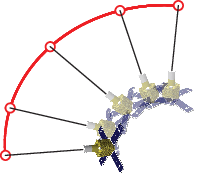}}
      \caption{\textsf{``Snapshots'' of the shortest-time (left) and industry-heritage (right) maneuvers reconstructed from TRACE's telemetry data of the 2010 flight demonstrations.}}
      \label{fig:TRACE_strobe}
   \end{figure}

\begin{figure}[h!]
      \centering
      {\includegraphics[scale=0.55]{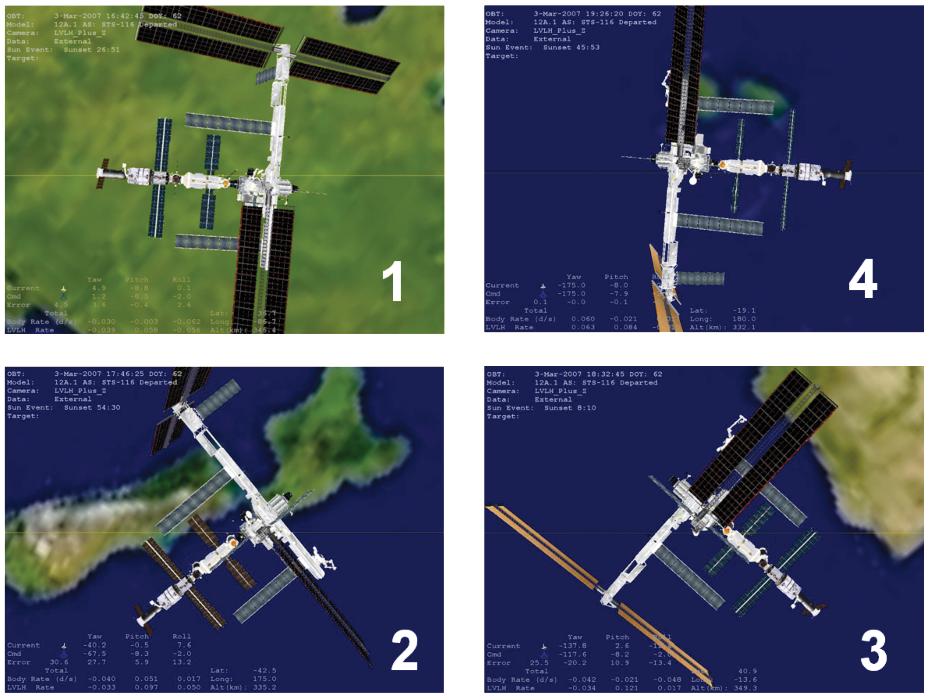}}
      \caption{\textsf{``Snapshots'' of Bedrossian's zero-propellant maneuver performed onboard the International Space Station; the visuals are reconstructions of the telemetry data of the 2006 flight of PS optimal control.}}
      \label{fig:ZPM_snapshots}
   \end{figure}

On August 10, 2010, a NASA space telescope called TRACE
executed the first ever minimum-time rotational maneuver performed in orbit~\cite{jgcd:trace-1,aas:trace-2,NASA_TRACE_article}.  This maneuver is the space analog of the classic Brachistochrone problem in which the spacecraft executes a longer path than an eigenaxis maneuver but reaches the goal faster~\cite{jgcd:rigid}; see Fig.~\ref{fig:STM_Brac_q3q2}.
%
%
The flight demonstration marked yet another first for NASA and another milestone for pseudospectral (PS) optimal control theory~\cite{SIAMnews, naz:firstISS}. ``Snapshots'' of the shortest-time maneuver, reconstructed from TRACE's telemetry data, are shown in Fig.~\ref{fig:TRACE_strobe}
%
%
alongside the industry heritage maneuver.   The TRACE flight demonstration was not the first time PS optimal control theory was used on orbit. In fact, PS optimal control theory debuted in flight on November 5, 2006 when NASA used it to implement Bedrossian's zero-propellant maneuver onboard the International Space Station~\cite{zpm:second,Bedrossian_2009}. ``Snapshots'' of this historic debut flight are shown in Fig.~\ref{fig:ZPM_snapshots}.
%
%

It is important to note that the zero-propellant maneuver was discovered, designed and implemented in orbit -- all using pseudospectral optimal control.  When this fact is juxtaposed with the fact that the flight implementation of the maneuver was performed on a 100 billion-dollar asset~\cite{GAO} with an international crew onboard, supreme confidence of technical success is essential.  That PS optimal control passed this high threshold for a NASA flight is an indication of how far the theory has evolved as a space technology.
%
%
Today, PS solutions are used to control everything from large industrial mechanisms~\cite{CDC:session} to quantum-mechanical systems~\cite{Li:pnas}.

In this paper, we review the theoretical foundations of PS optimal control that are critical for a successful practical implementation. Practical implementations, particularly manned flights, require mathematical guarantees that can be assured by numerical simulations and ground tests.  A good starting point for understanding the mathematical guarantees of PS optimal control is the notion of a Sobolev space, denoted as
$W^{m,p}$, as the space of all functions, $y$, whose $j$-th distributional
derivative, is in $L^p$ for all $0 \le j \le m$.
The Sobolev norm of $y$ is defined as,
\begin{equation}\label{eq:SobolevNorm}
\norm{y}_{W^{m,p}} :=\sum_{j=0}^{m} \norm{y^{(j)}}_{L^p} \nn
\end{equation}
where $\norm{y}_{L^p}$ is the usual Lebesgue norm,
\begin{equation}\label{eq:LpNorm}
\norm{y}_{L^p} := \left(\ds \int_\Omega \abs{y(t)}^p dt \right)^{1/p}\nn
\end{equation}
where $\abs{\cdot}$ denotes the absolute value.  
For notational
ease, we suppress the dependence of $W^{m,p}\ $ on vector-valued
functions as well as the domain $\Omega$. Throughout this paper, we will confine our discussions to state trajectories, $\bx(\cdot)$, that are in $W^{m_x,\infty},\ m_x \geq 1$ and control trajectories, $\bu(\cdot)$ that are in $W^{m_u,\infty},\ m_u \geq 0$.

\section{Practical Optimal Control: Challenges and Curses}
In order to set the stage, we review some of the basics of formulating and analyzing practical optimal control problems.  These basics are also a relevant starting point for pseudospectral optimal control theory.

\subsection{Search Space}

Let $\X \subset \real{N_x}$ and $\U \subset \real{N_u}$ be bounded sets that constrain the values of the state and control trajectories respectively, where $N_x, N_u \in \mathbb{N}$. In typical engineering optimal control problems~\cite{lu,rea,stanton,harada:1, hawkins,williams:terrain,RobSteve08} these sets are quite easily determined by the practical bounds imposed by the physics of the problem.  
For instance, for low-Earth-orbiting spacecraft, we can safely set the maximum distance to be the radius of the solar system. As discussed in Section 4.3, an assumption of such bounds allows us to formulate convergence results for PS optimal control. In a computational environment, we typically set $\X$ and $\U$ to be box constraints~\cite{DIDO_2007},
\begin{equation}
\X = \set{\bx \in \real{N_x}:\ \bx^L \le \bx \le \bx^U } \quad
\U = \set{\bu \in \real{N_u}:\ \bu^L \le \bu \le \bu^U }
\end{equation}
which are usually inactive unless the problem formulation specifically calls for such a constraint.

\subsection{Problem Formulation}
A practical optimal control problem may be posed as follows: Determine the state-control function pair, $t \mapsto (\bx \in \X,\ \bu \in \U)$, and clock times $\set{t_0, t_f}$ that solve the following problem,
\begin{eqnarray*}
&(\textsf{P}) \left\{
\begin{array}{lll}
\textsf{Minimize }  & {J}[\bx(\cdot), \bu(\cdot), t_0, t_f]\\
& \quad = E(\bx_0, \bx_f, t_0, t_f)
 + \ds \int^{t_f}_{t_0} {F}(\bx (t), \bu (t), t)\, dt \\
\textsf{Subject to} & \dot{\bx}(t)={\bf f}(\bx(t),\bu(t), t)\\
&\be^L \le \be\big(\bx_0, \bx_f, t_0, t_f \big) \leq \be^U \\
&\bh^L \le \bh(\bx(t), \bu(t), t) \le \bh^U
\end{array} \right.&
\end{eqnarray*}
where, $\bx_0 \equiv \bx(t_0)$, $\bx_f \equiv \bx(t_f)$, $\be^L, \be^U \in \real{N_e}$ are the lower and upper bounds on the endpoint function, $\be:\real{N_x} \times \real{N_x} \times \Real \times \Real \to \real{N_e} $, and $\bh^L, \bh^U \in \real{N_h}$ are the lower and upper bounds on the path function, $\bh:\real{N_x} \times \real{N_u} \times \Real \to \real{N_h} $. It is
assumed that $\be$, $\bh$ as well as all the other nonlinear functions, $F: \real{N_x} \times \real{N_u} \times \Real \to \Real$, $E:\real{N_x} \times \real{N_x} \times \Real \times \Real \to \Real$, and ${\bf f}: \real{N_x} \times \real{N_u} \times \Real \to \real{N_x}$, are continuously differentiable with respect to their arguments and that their gradients are Lipschitz-continuous over the
domain, $\X \times \U \times [t_0, t_f]$. We also assume that an optimal state-control function pair, $\set{\bx^*(\cdot),
\bu^*(\cdot)}$, exists in the Sobolev space  $W^{m_x,\infty}\times
W^{m_u,\infty}$ for some $m_x\geq 1$ and $m_u \ge 0$.


\subsection{Necessary and Sufficient Conditions}
The state dependent control space $\U(\bx, t)$ is defined by,
$$\U(\bx, t) := \set{\bu(t) \in \U:\ \bh^L \le \bh(\bx(t), \bu(t), t) \le \bh^U, \
\bx \in \X,\ t \in [t_0, t_f] }   $$
Thus, the lower Hamiltonian~\cite{vinter,Ross_BOOK}  is given by,
\begin{equation}\label{eq:hmc}
\mathcal{H}(\blam, \bx, t) := \min_{\bu \in \U(\bx,t)} H(\blam, \bx, \bu, t)
\end{equation}
where, $H$, defined by,
$$ H(\blam, \bx, \bu, t) := F(\bx, \bu, t) + \blam^T{\bf f}(\bx, \bu, t) $$
is the control Hamiltonian, and $\blam$ may be viewed as a parameter in $\real{N_x}$ that takes the role of the adjoint covector in Pontryagin's Principle or the (generalized) gradient of the Value Function, $V(\bx, t)$, in Bellman's Principle.  Thus, the sufficient condition for optimality is compactly given by the Hamilton-Jacobi-Bellman equation,
\begin{eqnarray*}
\mathcal{H}(\partial_xV(\bx, t), \bx, t) + \partial_tV(\bx,t) = 0 \qquad for \quad \ (\bx, t) \in \Omega
\end{eqnarray*}
which must be satisfied with appropriate boundary conditions.  Deferring a discussion of the boundary conditions to Section 2.4, note also that $\Omega$, the domain of $V$ is complicated by the presence of state constraints in the definition of Problem $P$.

The necessary conditions for Problem $P$ are given by Pontryagin's Principle which asserts the existence of covector functions, $t \mapsto (\blam, \bmu)$, and $\bnu$, that satisfy the following conditions\cite{Ross_BOOK},
\begin{eqnarray*}
& (\textsf{$P^\lambda$}) \left\{
\begin{array}{lrl}
& \dot{\bx}(t)-\partial_\lambda \overline{H}(\bmu(t), \blam(t), \bx(t),\bu(t), t) = &0\\
& \dot{\blam}(t)+\partial_x \overline{H}(\bmu(t), \blam(t), \bx(t),\bu(t), t) = &0\\
&\bh^L \le \bh(\bx(t), \bu(t), t) \le & \bh^U\\
& \partial_u \overline{H}(\bmu(t), \blam(t), \bx(t),\bu(t), t) = &0\\
&\be^L \le \be\big(\bx_0, \bx_f, t_0, t_f \big) \leq &\be^U \\
&\blam(t_0) + \partial_{x_0}\overline{E}(\bnu,\bx_0,\bx_f, t_0, t_f) = &0\\
&\blam(t_f) - \partial_{x_f}\overline{E}(\bnu,\bx_0,\bx_f, t_0, t_f) = &0\\
& \overline{H}[@ t_0]- \partial_{t_0}\overline{E}(\bnu,\bx_0,\bx_f, t_0, t_f) = &0 \\
& \overline{H}[@ t_f]+ \partial_{t_f}\overline{E}(\bnu,\bx_0,\bx_f, t_0, t_f) = &0 \\
&\bmu\ \bot & \bh\\
&\bnu\ \bot &\be
\end{array} \right. &
\end{eqnarray*}
where $\overline{H}$ is the Lagrangian of the Hamiltonian,
$$ \overline{H}(\bmu, \blam, \bx, \bu, t) := H(\blam, \bx, \bu, t) + \bmu^T\bh(\bx, \bu, t)   $$
$\overline{H}[@ t_0]$ is a shorthand notation for $\overline{H}(\bmu(t_0), \blam(t_0), \bx(t_0),\bu(t_0), t_0) $ with $\overline{H}[@ t_f]$ defined similarly,
$\overline{E}$ is the endpoint Lagrangian,
$$ \overline{E}(\bnu, \bx_0, \bx_f, t_0, t_f) := E(\bx_0, \bx_f, t_0, t_f) + \bnu^T\be(\bx_0, \bx_f, t_0, t_f)  $$
and $\bmu \bot \bh$ is a shorthand notation for the complementarity condition\cite{Ross_BOOK},
$$\mu_i \left\{\begin{array}{ccc}
               \le 0            & if & h_i(\bx(t), \bu(t), t) = h_i^L \\
               =  0             & if & h_i^L < h_i(\bx(t), \bu(t), t) < h_i^U \\
               \ge 0            & if & h_i(\bx(t), \bu(t), t) = h_i^U \\
               unrestricted     & if & h_i^L = h_i^U
             \end{array}
   \right.
$$
with $\bnu \bot \be$ defined similarly.

\subsection{Challenges and Curses}
Problem $P$ is a standard formulation of practical engineering problems; see~\cite{lu,rea,stanton,harada:1, hawkins,williams:terrain,RobSteve08} and Section 5 of this paper. It is clear that even assembling the sufficient conditions for optimality is challenging because it involves the following steps:
\begin{enumerate}
\item A closed form expression for $\mathcal{H}$, which requires a solution to (\ref{eq:hmc}), a formidable task in and of itself;
\item A functional description of $\Omega$ which is a subset of $\real{N_x}\times \Real$, a challenging task in the presence of state constraints; and,
\item A formulation of the appropriate boundary conditions for $V(\bx, t)$, made difficult by the nature of mixed boundary conditions and inequalities.
\end{enumerate}
These practical difficulties are not the topic of typical mathematical research because even if all these challenges are met, there is the overhanging doom of the curse of dimensionality coupled with the issues related to the nondifferentiability of $V$~\cite{vinter}.

In sharp contrast to the sufficient conditions, the necessary conditions offer a more tractable approach, but is not without its own challenges. The necessary conditions are tractable, at least in their assemblage; see Problem $P^\lambda$. Solving for the necessary conditions is motivated by the possibility that the extremal is an optimal solution.  Problem $P^\lambda$ is a generalized root-finding problem that can be characterized as a boundary value problem (BVP) with the following challenges:
\begin{enumerate}
\item $P^\lambda$ is a system of differential-algebraic equations (DAEs). DAEs have significantly more theoretical and computational challenges than a system of ordinary differential equations.

\item Solving a BVP is fundamentally a more difficult task than an initial value problem (IVP). In fact, a BVP may not have a solution whereas a corresponding IVP may have one or more solutions. $P^\lambda$ is not only a BVP, but its solution is composed of switches in $t \mapsto (\bu, \blam, \bmu) $ with no \textit{a\,priori} information on where and how many switches the solution may have.
\end{enumerate}
Should these challenges be overcome, there remains the overarching curse of sensitivity~\cite{Ross_BOOK} arising from the symplectic structure of the Hamiltonian system.  That is, the inevitability of instability in the numerical solution of $P^\lambda$ --- a consequence of the eigenvalues (of the linearized system) lying on both sides of the imaginary axis~\cite{Ross_BOOK,RossGongBook,Kalman:difficulty}.  

\section{Pseudospectral Foundations}
Pseudospectral optimal control theory, a term coined by Ross, is a joint theoretical-computational framework for solving Problem $P$.  It is founded on the fact that the space of polynomials is dense over a compact interval, which is a direct consequence of Weierstrass' famous approximation theorem:
\begin{theorem}[Stone-Weierstrass]\label{thm:sw}
Let $\Real \supset [t_0, t_f] \mapsto y(t) \in \Real$ be a continuous function.  Then, there exists an algebraic polynomial sequence, $\set{t \mapsto y^N}_{N=0}^\infty$ on $[t_0, t_f]$ such that,
$$ \lim_{N \to \infty} \sup_{t \in [t_0, t_f]}\abs{y(t)-y^N(t)} = 0  $$
\end{theorem}
Because the state trajectory, $\bx(\cdot)$, is assumed to be in $W^{1,\infty}$, it follows that there exists a polynomial sequence, $\set{t \mapsto \bx^N}_{N=0}^{\infty}$ that converges to $\bx(\cdot)$ in $L^\infty$.
%
%
In PS theory, we aim to find  $\bx^N(t)$ by judiciously combining the following four elements~\cite{RossGongBook, Fahroo_2008a, elnagar:first, elnagar:cheb, fahroo:first}:
(\emph{i}) domain transformation, (\emph{ii}) interpolation, (\emph{iii}) differentiation and (\emph{iv}) integration.

\subsection{Domain Transformation}
For finite-horizon optimal control problems, we map the physical domain $[t_0, t_f] \ni t$ to a computational domain $[-1, 1] \ni \tau$ by means of an affine transformation,
\begin{equation}\label{eq:time-finite}
t = \Gamma(\tau) =\left(\frac{t_f+t_0}{2}\right) + \left(\frac{t_f-t_0}{2}\right)\,\tau \quad  \Leftrightarrow \quad   \tau = \Gamma^{-1}(t) = \left(\frac{2}{t_f-t_0}\right) t - \left(\frac{t_f+t_0}{t_f - t_0}\right)
\end{equation}
Other domain transformations techniques~\cite{Fahroo_2008a,adaptive, AAS:gauss} can be used for enhanced computational efficiency.

For infinite-horizon optimal control problems, we map the (semi-) infinite domain, $[t_0, \infty)$, to the half open interval, $[-1, 1)$, using the transformation,
\begin{equation}\label{eq:time-inf}
t=\Gamma(\tau) = \frac{1+\tau}{1-\tau} \ \Leftrightarrow \
\tau = \Gamma^{-1}(t) = \frac{t-1}{t+1}
\end{equation}
As in the finite horizon case, other domain transformational techniques can be used as well~\cite{Fahroo_2008a, AAS:gauss, Radau-GNC05, acc:stability, PS:Radau}.

To maintain notional clarity, we abuse notation in the rest of this paper by using the same symbol $t$ to mean $ t \in [t_0, t_f]$ as well as $t \in [-1, 1]$. The domain to which we are referring will be apparent from the context.

\subsection{Interpolation}

\begin{figure}[h!]
      \centering
      {\includegraphics[scale = 0.5]{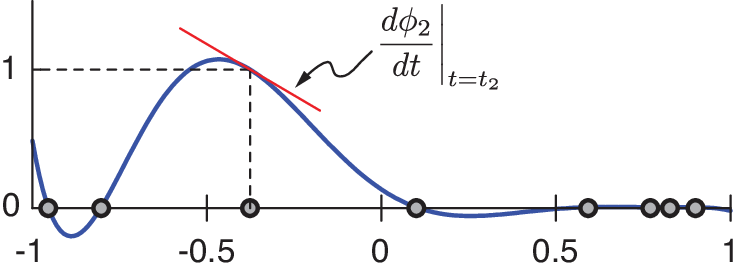}}
      \caption{\textsf{An arbitrary grid, $\pi^N$, for $N = 7$, and the Lagrange polynomial $\phi_2(t)$ over $\pi^N$.}}
      \label{fig:Arbitrary_grid}
   \end{figure}

\begin{figure}[h!]
      \centering
      {\includegraphics[scale = 0.5]{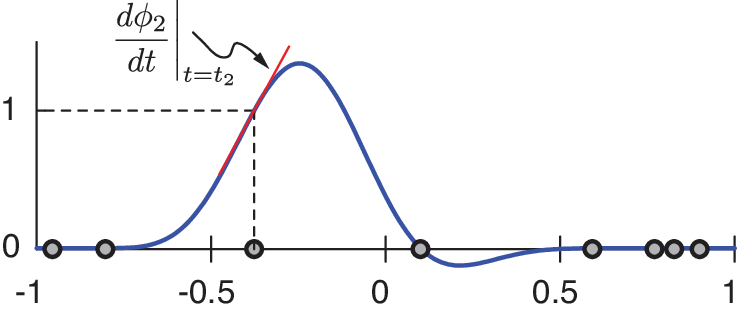}}
      \caption{\textsf{A weighted interpolating polynomial over the same arbitrary grid, $\pi^N$ of Fig.~\ref{fig:Arbitrary_grid}.}}
      \label{fig:Weighted_Poly}
   \end{figure}

Let $\pi^N := \set{t_j, j=0,...,N}$ be a set of arbitrary but distinct points or ``nodes''
on the interval $[-1, 1]$.
A standard Lagrange interpolant  can be written as,
\begin{equation}\label{eq:interp}
y^N (t)=\sum_{j=0}^{N} \phi_j
(t)y_j,\qquad -1 \leq t \leq 1
\end{equation}
where $\phi_j(t)$ is the $N^{th}$-order
Lagrange interpolating polynomial (see Fig.~\ref{fig:Arbitrary_grid}),
\begin{equation}\label{eq:Lpoly-basic}
\phi_j (t)=\prod_{i=0, i \ne j}^{N}\frac{t-t_i}{t_j-t_i}.
\end{equation}
%
%
%
It is apparent that $\phi_j(t)$ satisfies the Kronecker relationship
$\phi_j(t_k)=\delta_{jk}.$ This implies that
\begin{equation}\label{eq:interp-pts}
y_j=y^N(t_j), \,\,\, j=0,...,N.
\end{equation}
In pseudospectral optimal control, we use weighted interpolants of the form,
\begin{equation}\label{eq:interpW}
y^N (t)=\sum_{j=0}^{N} \frac{W(t)}{W(t_j)}\phi_j
(t)y_j,\qquad -1 \leq t \leq 1
\end{equation}
where  $W$ is a positive
weight function. The simplest weight function is $W(t) \equiv 1$ in which case (\ref{eq:interpW}) simplifies to (\ref{eq:interp}).
The main difference between (\ref{eq:interpW}) and (\ref{eq:interp}) is that the weight function provides a way to shape the derivative of the interpolants while maintaining the Kronecker property (see Fig.~\ref{fig:Weighted_Poly}).
%
%

\subsection{Differentiation}
It is straightforward to show that,
\begin{equation}\label{eq:poly-diff}
\left.\frac{d y^N(t)}{dt}\right|_{t_i}=\sum_{j=0}^{N}
\frac{y_{j}}{W(t_j)}[\dot W(t_i)\delta_{ij} +W(t_i)\widehat{D}_{ij}] =
\sum_{j=0}^{N} D_{ij} y_j
\end{equation}
where $D_{ij}$ is a shorthand notation for the
$W$-weighted differentiation matrix,
\begin{equation}\label{eq:D-def-Wt}
D_{ij} := \frac{[\dot W(t_i)\delta_{ij}
+W(t_i)\widehat{D}_{ij}]}{W(t_j)}
\end{equation}
and $\widehat{D}_{ij}$ is unweighted differentiation matrix given by,
\begin{equation}\label{eq:D-def-unWt}
\widehat{D}_{ij} :=\left.\frac{d\phi_j(t)}{dt}\right|_{t=t_i}
\end{equation}
Thus, when $W(t) \equiv 1$, we have
\begin{equation}\label{eq:D-def-1}
D_{ij} = \widehat{D}_{ij}
\end{equation}

\subsection{Integration}
Similar to differentiation, we have,
$$\int_{-1}^{1}y^N(t)\, dt =\sum_{j=0}^{N} \int_{-1}^{1}\frac{W(t)}{W(t_j)}\phi_j
(t)\, dt \ y_j $$
Thus, we get,
\begin{equation}\label{eq:quad-good}
\int_{-1}^{1}y^N(t)\, dt := \sum_{j=0}^N  w_{j} y_{j}
\end{equation}
where $w_j$ are the integration weights given by
\bea
w_j & = & \int_{-1}^1\frac{W(t)}{W(t_j)} \phi_j(t)\, dt\qquad j=0,1,...,N \label{eq:weights}
\eea

\begin{remark}\label{rm:summary}
A pseudospectral method is completely determined by just three design parameters, $\Gamma$, $\pi^N$, and $W$. Standard choices for these parameters are~\cite{RossGongBook, Fahroo_2008a,PS:Radau}:
\begin{enumerate}
\item Linear and bi-linear for $\Gamma$; linear for finite-horizon and bi-linear for infinite-horizon optimal control problems;
\item Gauss-Lobatto and Gauss-Radau for $\pi^N$; Gauss-Lobatto for finite-horizon and Gauss-Radau for infinite-horizon optimal control problems; and,
\item $W(t) \equiv 1$ (Gauss-Lobatto) and $W(t) = 1- t$ (Gauss-Radau).
\end{enumerate}
\end{remark}

%

%
%
%

\section{Overview of Pseudospectral Optimal Control}

While Problem $P$ is a very practical optimal control problem, 
the development of pseudospectral theory for Problem $P$ generates a number of bookkeeping issues that are related to the various dimensions of the decision variables and constraints (such as $N_x, N_e $ etc.).  Additional well-known bookkeeping issues arise due to the use of lower and upper bounds to generalize equality and inequality constraints. As a means to prevent these distractions and limit the discussions to the core issues of pseudospectral optimal control theory, we consider the following representative problem,
\begin{eqnarray*}
&x \in \X \subset \Real, \quad u \in \U \subset \Real   &\\
& (\textsf{$B$}) \left\{
\begin{array}{lrl}
\textsf{Minimize } & J[x(\cdot), u(\cdot)] =& E(x(-1),x(1)) \\
\textsf{Subject to}& \dot x(t) =& f(x(t), u(t))  \\
& e(x(-1), x(1))  = & 0
\end{array} \right. &
\end{eqnarray*}
All critical concepts of PS theory discussed for Problem $B$ easily extend to Problem $P$.
Following the discussions of Section 2.3, Problem $B$ generates the dualized Problem $B^{\lambda}$ defined by,
%
\begin{eqnarray*}
& (\textsf{$B^\lambda$}) \left\{
\begin{array}{lrl}
\textsf{Find} & t
\mapsto ( x,u,\lambda) \textsf{ and } \nu \\
\textsf{Such that} & \dot x(t) - f(x(t), u(t)) =&0 \\
& e(x(-1), x(1))  = & 0 \\
&\dot{\lambda}(t) + {\partial_x {H}}(\lambda(t),x(t),u(t)) = & 0\\
&{\partial_u {H }}(\lambda(t),x(t),u(t)) =&0\\
&\lambda(-1) + \ds\frac{\partial {\overline{E}}}{\partial x(-1)}(\nu,x(-1),x(1)) = &0\\
&\lambda(1) - \ds\frac{\partial {\overline{E}}}{\partial
x(1)}(\nu,x(-1),x(1)) = &0
\end{array} \right. &
\end{eqnarray*}

\subsection{Theoretical Overview}


Pseudospectral optimal control theory is founded on generating a convergent sequence of interpolating functions that satisfy the covector mapping principle of Ross and Fahroo~\cite{prague,lncis,acc:hybrid,cmp:history,cmp:cdc2006}.

The sequence generator for $t \mapsto x^N $ is given by,
\begin{equation}\label{eq:xn}
x^N(t)=\sum_{j=0}^{N} \frac{W(t)}{W(t_j)} x_j \phi_j(t)
\end{equation}
where $W$ is an arbitrary weight function that determines a particular design of a pseudospectral method.  The sequence generator for $t \mapsto u^N $ is given by,
\begin{equation}
u^N(t)=\sum_{j=0}^{N} u_j \psi_j(t)
\end{equation}
where $\psi_j(t)$ is a special interpolating function whose property will be apparent shortly. The coefficients $x_i, u_i$ of the interpolating functions are solutions to Problem $B^N$ defined by~\cite{Fahroo_2008a,TAC:linearizable},
\begin{eqnarray*}
&x_i \in \X, \quad u_i \in \U   &\\
& (\textsf{$B^N$}) \left\{
\begin{array}{lrl}
\textsf{Minimize } & J^N [x(\pi^N), u(\pi^N)] &=
 E(x_0, x_N)\\ 
\textsf{Subject to}& \norm{\ds\sum_{j=0}^{N}D_{ij}x_j -f(x_i, u_i)} &\le  \delta^N \quad i = 0, 1, \ldots, N  \\
& \norm{e(x_0, x_N)}  &\le \delta^N
\end{array} \right. &
\end{eqnarray*}
where $x(\pi^N) =(x_0, x_1, \ldots, x_N) \in \real{N+1} $, $u(\pi^N)$ is defined similarly, $\pi^N$ is an arbitrary grid, and $\delta^N$ is a positive number that depends on $N$ with the property that $\delta^N \to 0$ as $N \to \infty$.
\begin{remark}
The control function, $t \mapsto u^N$, is
not necessarily a Lagrange interpolant;
rather, $u^N(t)$ is a generalized interpolant that satisfies,
\begin{equation}
\norm{\dot x^N(t) -f(x^N(t), u^N(t))} \rightarrow 0 \quad  as \quad \delta^N \rightarrow 0
\end{equation}
This aspect of PS control is exploited in the Bellman pseudospectral method~\cite{bellmanPaper:JGCD,bellmanPaper:aas} where information in between the nodes is obtained by a first-principles application of the principle of optimality.
\end{remark}
The covector mapping principle (CMP) was the first major result developed in PS optimal control theory~\cite{RossGongBook, fahroo:first,prague, lncis, cmp:history, fahroo:LGL}. Early results~\cite{fahroo:first,fahroo:LGL} were based on convergence assumptions; thanks to convergence theorems developed by Gong et al.~\cite{TAC:linearizable, coap:cmp} and Kang et al.~\cite{kang:disc, krg:convergence}, the CMP is now firmly established as a fundamental result in PS theory. This result can be articulated as follows:\\
%

\noindent\textit{\textbf{Covector Mapping Principle}} \textbf{(Ross-Fahroo)} \
\textit{Let the sequence of function pairs $t \mapsto \{x^N, u^N\}_{N=N_0}^\infty$ converge to a solution of Problem $B$. Then, there exist multipliers for Problem $B^N$ that map to the coefficients of some interpolating functions that converge to covector functions that solve Problem $B^\lambda$}.\\
%
%

Let the interpolating covector function  $t \mapsto \lambda^N$ be given by,
\begin{equation}\label{eq:ln}
\lambda^N(t)=\sum_{j=0}^{N} \frac{W^*(t)}{W^*(t_j)} \lambda_j
\phi_j(t)
\end{equation}
where $W^*$ is some (dual) weight function related to a particular design of a pseudospectral method~\cite{Fahroo_2008a}.  The coefficients $\lambda_i$ are solutions to Problem $B^{\lambda N}$ defined by~\cite{AAS:gauss,coap:cmp},
\begin{eqnarray*}
& (\textsf{$B^{\lambda N}$}) \left\{
\begin{array}{lrl}
\textsf{Find} & x_i \in \X, u_i \in \U, \lambda_i \in \Real \textsf{ and } \nu^N & \quad\  \qquad\, i = 0, 1, \ldots, N\\
\textsf{Such that} &\norm{\ds\sum_{j=0}^{N}D_{ij}x_j - f(x_i, u_i)} & \le \delta^N \qquad i = 0, 1, \ldots, N  \\
& \norm{e(x_0, x_N)} &  \le \delta^N \\
&\norm{\ds\sum_{j=0}^{N}D^*_{ij}\lambda_j +{\partial_{x_i} {H}}(\lambda_i,x_i,u_i)} &\le \delta^N  \qquad i = 0, 1, \ldots, N  \\
&\norm{{\partial_{u_i} {H }}(\lambda_i,x_i,u_i)} & \le \delta^N \qquad i = 0, 1, \ldots, N  \\
&\norm{\lambda_0 + \partial_{x_0}\bar{E}(\nu,x_0,x_N) } & \le \delta^N\\
&\norm{\lambda_N - \partial_{x_N}\bar{E}(\nu,x_0,x_N) }
& \le \delta^N
\end{array} \right. &
\end{eqnarray*}
where $D^*_{ij}$ is the differentiation matrix associated with $W^*$.
Let $\tlam_i, i= 0, 1, \ldots, N$ and $\tnu^N$ be a set of multipliers for Problem $B^N$. Then, from the CMP, if the sequence $t \mapsto \{x^N, u^N\}_{N=N_0}^\infty$ converges, we can find maps, $\tlam_i \mapsto \lambda_i$ and $\tnu^N \mapsto \nu^N$, such that $\{ t \mapsto \lambda^N, \nu^N \}_{N=N_0}^\infty$ converges to a solution of Problem $B^\lambda$.

It is possible to obtain some general results for the case of a linear map.
\begin{proposition}
Let, $\widetilde{\Lambda} = [\widetilde{\lambda}_0, \widetilde{\lambda}_1, \cdots, \widetilde{\lambda}_N]^T$ and $\Lambda = [\lambda_0, \lambda_1, \cdots, \lambda_N]^T$. Let $P$ be an $(N+1) \times (N+1)$ matrix such that,
$$ \Lambda = P \widetilde{\Lambda}$$
If $\widetilde{\Lambda}$ and $\Lambda$ satisfy the CMP, and $\partial_x f(x(t), u(t), t)$ is not a constant, then $P$ is an invertible diagonal matrix.
\end{proposition}

\begin{remark}
Proposition 1 generalizes the results obtained for the ``big two'' techniques, namely the Legendre\cite{lncis} and Chebyshev\cite{fahroo:cheb-jgcd,cheb:cdc-mapping} pseudospectral methods.  In both cases, the diagonal elements of $P$ are the reciprocals of the integration weights, $w_i$. Thus, for the Legendre PS method, $w_i$ are the Gauss-Lobatto weights while for the Chebyshev PS method, $w_i$ are the Clenshaw-Curtis weights.
\end{remark}

\begin{figure}[h!]\centering\scalebox{0.5}
      {\includegraphics[scale = 0.5]{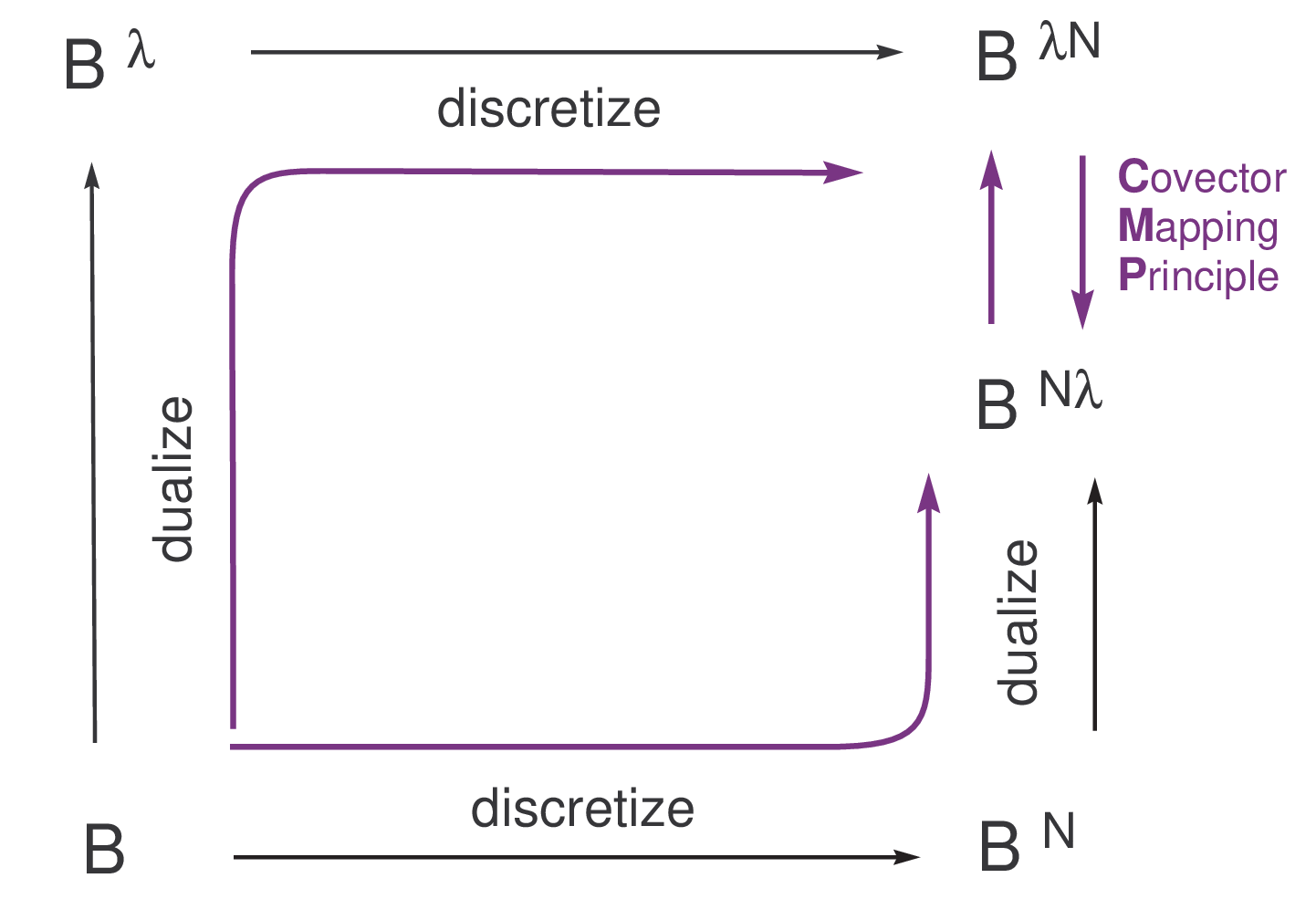}}
      \caption{\textsf{Illustrating the Covector Mapping Principle as a commutative diagram of Ross and Fahroo.
      }}
      \label{fig:CMP-commute}
   \end{figure}

An informal statement of the CMP is that it is possible to commute the operations of dualization and discretization.  This fact is illustrated by the commutative diagram shown in Fig.~\ref{fig:CMP-commute}.
%
%

\subsection{Design of Pseudospectral Methods}
The design of Problem $B^N$ is completely determined by $\delta^N$, $\pi^N$ and $W$.  The choices of these parameters determine the success or failure of a pseudospectral method.  Conditions for choosing these parameters are determined by convergence theorems outlined in Section 4.3. In this section, we illustrate some ``corners'' of the ``design space.''

\subsubsection{Choosing the consistency parameter, $\delta^N$}
The following theorem~\cite{TAC:linearizable,coap:cmp, arbGrid} provides guidelines on selecting $\delta^N$:
%
\begin{theorem}[Gong-Kang-Ross-Fahroo]  \label{T-exist}
Given any feasible solution, $t \mapsto (x,u)$, for Problem B, suppose
$x(\cdot) \in W^{m_x, \infty}$ with $ m_x \ge 2$. Then, there exists a positive
integer $N_0$ such that, for all $N>N_0$, Problem $B^N$ has feasible
solutions with $\delta^N = (N-1)^{3/2-m_x}$.
Furthermore, a feasible solution satisfies $u_i \ = \ u(t_i)$ and
\bea
\| x_i-x(t_i)\|_\infty & \leq & L(N-1)^{1-m_x},\nn
\eea
for all $i=0,\ldots,N$ where $L$ is a positive constant independent of $N$.
\end{theorem}

\subsubsection{Choosing the grid, $\pi^N$}
\begin{figure}[h!]\centering
      {\includegraphics[width=2.5in]{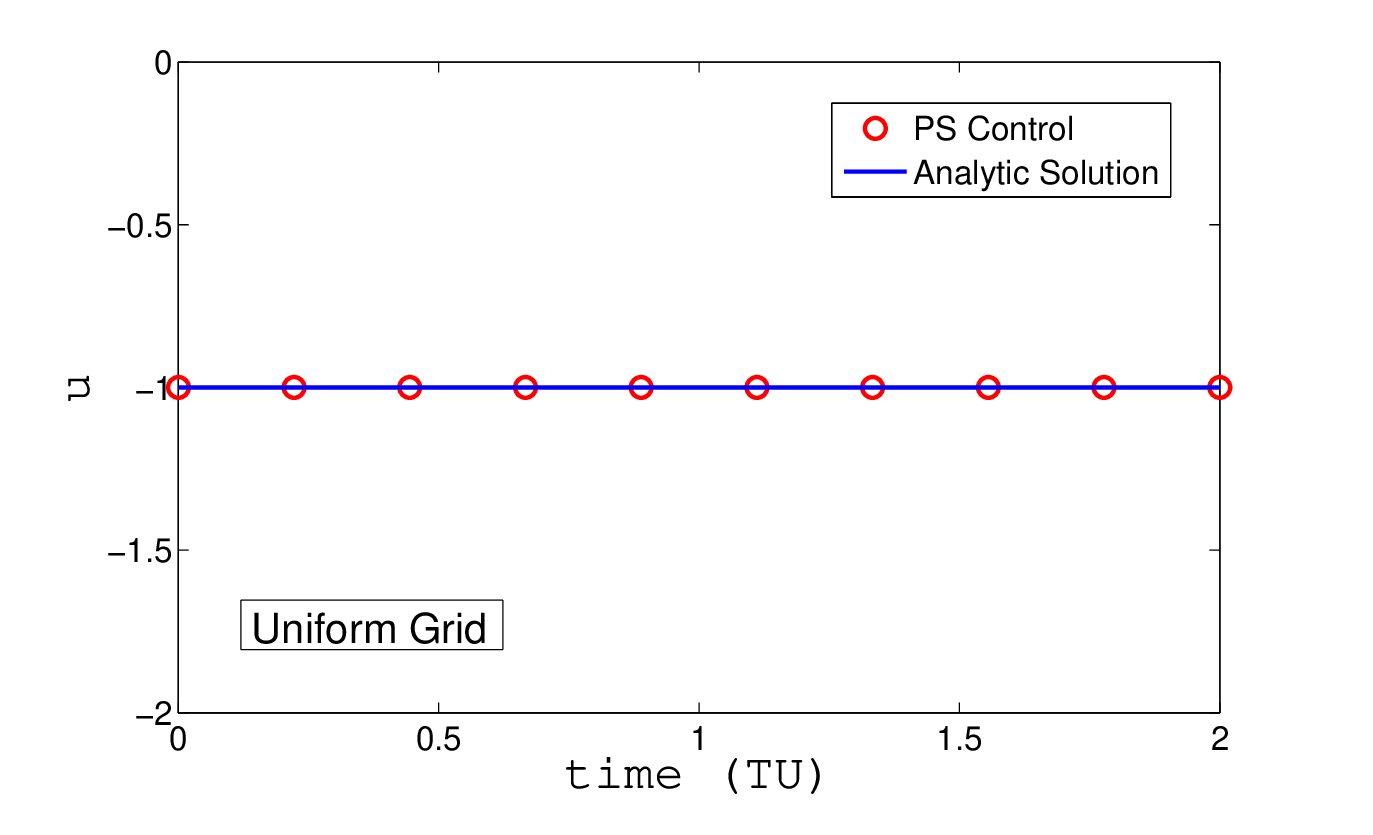}}
      \caption{\textsf{PS control solution for a uniform grid of 10 points for Problem $E_1$.}}
      \label{fig:Uniform10}
   \end{figure}

\begin{figure}[h!]\centering
      {\includegraphics[width=2.65in]{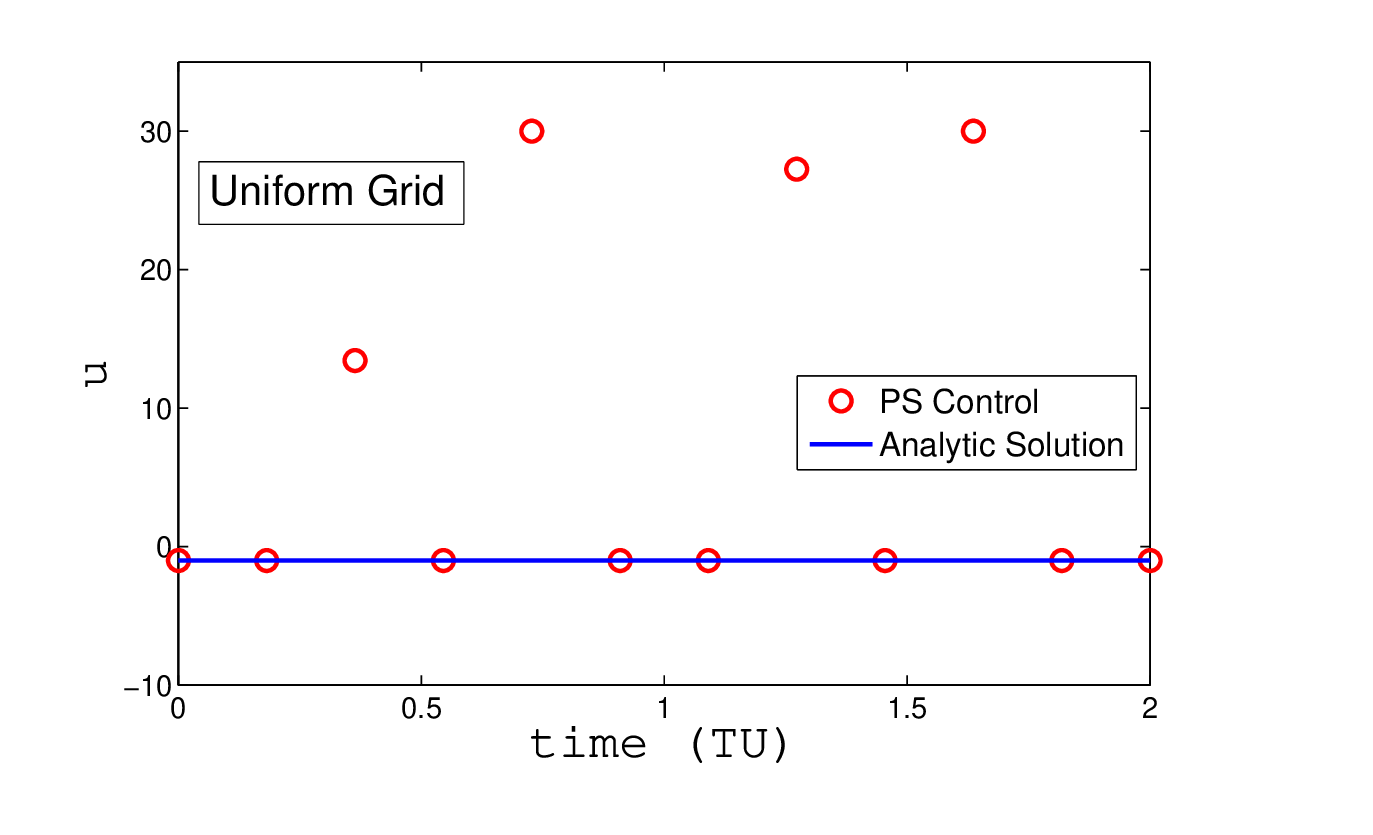}}
      \caption{\textsf{PS control solution for a uniform grid of 12 points for Problem $E_1$.}}
      \label{fig:Uniform12}
   \end{figure}

The simplest choice of a uniform grid, $\pi^N = \pi^N_{uniform}$, is the worst possible choice for $\pi^N$. To illustrate this point, consider the following example~\cite{arbGrid},
\begin{eqnarray*}
&x \in \Real, \quad u \in \U = \set{u \in \Real: u \ge -1}   &\\
& (\textsf{$E_1$}) \left\{
\begin{array}{lrl}
\textsf{Minimize } & J[x(\cdot), u(\cdot)] =& x(2) \\
\textsf{Subject to}& \dot x(t) =& u(t)  \\
& x(0)  = & 0
\end{array} \right. &
\end{eqnarray*}
It is straightforward to show that this problem has a unique closed-form control given by, $u^*(t) = -1$.
The PS control solution for $\pi^N_{uniform}$ for $N=10$ is shown in Fig.~\ref{fig:Uniform10}.
%
%
At first glance, the solution is quite impressive when compared to the analytic solution; however, when $N =12$, it is clear from Fig.~\ref{fig:Uniform12} that the solution diverges and hence does not converge as $N \to \infty$.
%
%

The reason the solution diverges is that the selection of $\pi^N$ must be limited to those grid points for which the integration weights, $w_i$, given by (\ref{eq:weights}) are positive.  This assumption is used in the proof of PS convergence theorems discussed in Section 4.3.  For a uniform grid, $w_i$ are positive for $N=10$ but when $N \ge 11$, at least one of the weights is negative.

Thus, it is clear that we need a more sophisticated approach for choosing $\pi^N$ than a simple uniform distribution.

\begin{figure}[h!]\centering
      {\includegraphics[width=2.5in]{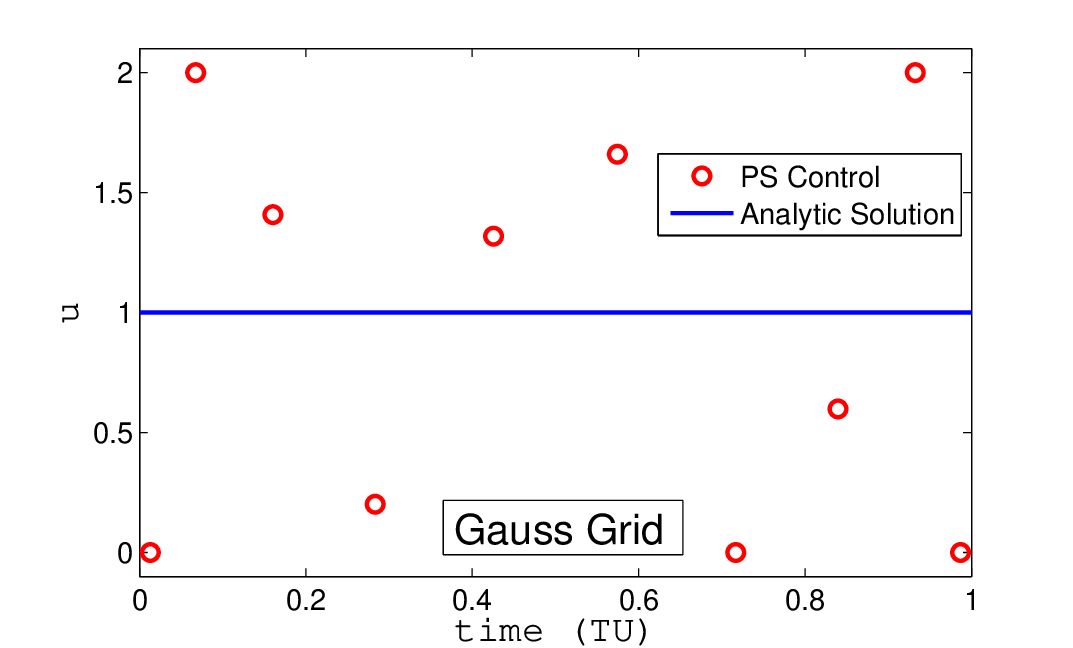}}
      {\includegraphics[width=2.5in]{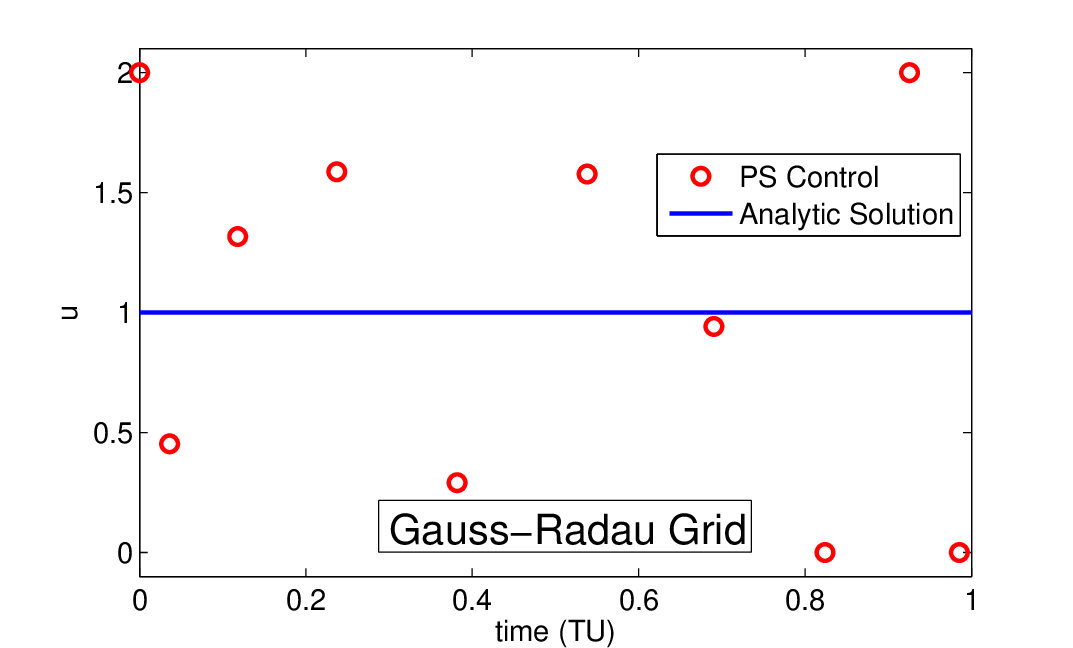}}
      {\includegraphics[width=2.5in]{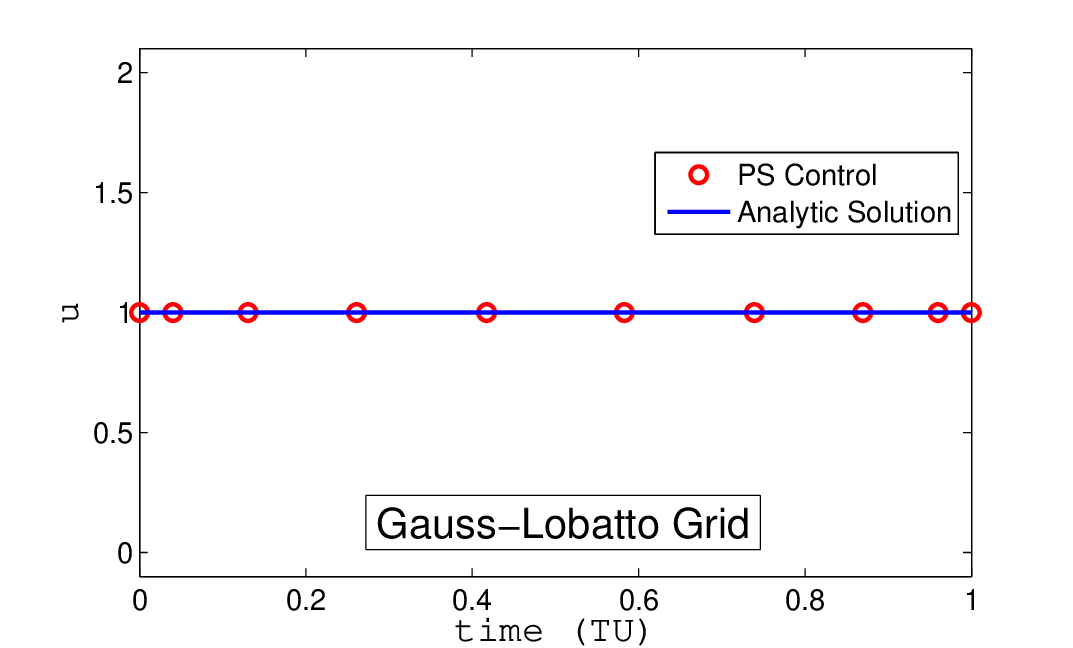}}
      \caption{\textsf{Divergence and convergence of PS control solutions for three Gaussian grids of 10 points for Problem $E_2$.}}
      \label{fig:Gaussian10}
   \end{figure}

%
\begin{proposition}\label{prop:gauss}
If $\pi^N$ is Gaussian, and $W(t) \equiv 1$, all $w_i, i = 0, \ldots, N$ are positive.
\end{proposition}
Based on this result, it is quite tempting to choose $\pi^N$ to be the simplest of the Gaussian grids, namely Gauss. The following example~\cite{jgcd:tn-no-conv} reveals that the design space within the Gaussian grids must be further constrained to Gauss-Lobatto points:
\begin{eqnarray*}
&\bx \in \X = \set{(x_1, x_2) \in \real 2: x_2 \ge 0 }, \quad u \in \U = \set{u \in \Real: 0 \le u \le 2 }    &\\
& (\textsf{$E_2$}) \left\{\begin{array}{lrl} \textsf{Minimize } &
J[\bx(\cdot), u(\cdot)]
=&\displaystyle \int_0^1 x_2(t) u(t) \ dt  \\
\textsf{Subject to} & \dot x_1(t) = & x_2(t)\\
& \dot x_2(t) = & - x_2(t)+u(t)\\
&(x_1(0),x_2(0))= & (0,1) \\
&(x_1(1),x_2(1))= & (1,1)
\end{array}\right.&
\end{eqnarray*}
The closed-form solution to Problem $E_2$ is given by,
\begin{eqnarray*}
x_1^*(t) = t, \qquad
x_2^*(t) = 1, \qquad
u^*(t) = 1
\end{eqnarray*}
PS control solutions for $N=10$ for three Gaussian grids (Gauss, Gauss-Radau and Gauss-Lobatto) are shown in Fig.~\ref{fig:Gaussian10}.
%
%
\textit{It is clear that only the solution over a Gauss-Lobatto grid converges}. This is because choosing $\pi^N$ to be Gaussian is not a sufficient condition for convergence; the grid must be paired correctly with a proper interpolation weight function, $W$.


\subsubsection{Choosing the interpolation weight function, $W$}
In a series of papers~\cite{Fahroo_2008a,AAS:gauss,Radau-GNC05,PS:Radau} Fahroo and Ross showed that a proper pairing of $\pi^N$ and $W$ are essential for convergence. They showed that these pairings are intimately connected to the type of boundary conditions stipulated by the problem, and that the uniform weight function, $W(t) \equiv 1$, cannot be associated with Gauss or Gauss-Radau grids; however, if one chooses nonuniform weight functions, such as $W(t) = 1- t$ or $W(t) = 1 - t^2$, then $\pi^N$ may be chosen to be Gaussian with appropriate caveats. These results are summarized in Table 1.


\begin{table}[h!]
\centering
\begin{tabular}{l c c r}
  \hline\hline
  W & $\pi^N$ & Boundary Conditions & Typical Horizon \\
  \hline
  1         & Gauss-Lobatto & Any           & Finite \\
  $1-t$     & Gauss-Radau   & Fixed-Free    & Infinite \\
  $1-t^2$   & Gauss         & Free-Free     & NA \\
  \hline
\end{tabular}
\caption{Sample weight functions and their applicability to PS optimal control.}\label{table1}
\end{table}

Note that Table 1 also associates a horizon for the $(\pi^N, W)$ pair. For the Gauss-Radau grid, the bilinear map (see Section 3 and Remark 1) or its equivalent is assumed for the domain transformation.




\subsection{Overview of Existence, Convergence and Consistency Theorems}
Theorems for pseudospectral optimal control are based on a key foundation that the finite-dimensional state and control spaces, $\X$ and $\U$, are bounded.  This is required to ensure the existence of the limits of sequences. This aspect of the foundation maps very well to the practical nature of the original problem formulation (see Problem $P$ in Section 2.1) which requires $\bx(t) \in \X$ and $\bu(t) \in \U$ for all $t$.  In addition, this formulation serves the computational need for limiting the search space in optimization.

A theoretical foundation for the existence of a solution for Problem $B^N$ comes from Weierstrass' Existence Theorem,
\begin{theorem}[Weierstrass]\label{thm:Weierstrass}
Let $\mathbb{Y} \subset \real{N_y}$ be a nonempty compact set and $C: \mathbb{Y} \to \Real$ be a continuous function; then, there exists a $y^* \in \mathbb{Y}$ such that $C(y^*) \le C(y)$ for all $y \in \mathbb{Y}$.
\end{theorem}
From Theorem \ref{T-exist}, it can be argued~\cite{TAC:linearizable,coap:cmp,kang:disc,arbGrid,cheb:cdc-mapping} that the feasible set for Problem $B^N$ is nonempty and compact; hence, for any given $N$, an optimal solution $ t \mapsto (x^{*N}, u^{*N})$ exists because $J^N$ is continuous (over a compact domain).

The theoretical basis for convergence is the Arzel\`{a}-Ascoli theorem,
\begin{theorem}[Arzel\`{a}-Ascoli]\label{thm:AA}
Let $\set{t \mapsto y^N}_{N=0}^\infty$ be a sequence of continuous functions defined over the closed interval $[-1, 1]$. If this sequence is uniformly bounded and uniformly equicontinuous, then it admits a subsequence which converges uniformly.
\end{theorem}
Thus, if the conditions of the Arzel\`{a}-Ascoli theorem are satisfied, we can select a subsequence, $\set{t \mapsto (x^{*n}, u^{*n})}_{n=N_0}^\infty$, of $\set{t \mapsto (x^{*N}, u^{*N})}_{N=N_0}^\infty$ such that the following limits
\begin{eqnarray}
\lim_{n \to \infty} x^{*n}(t)  = x^\infty(t)\qquad
\lim_{n \to \infty} u^{*n}(t)  = u^\infty(t)
\end{eqnarray}
converge uniformly. For the Legendre PS approach, the conditions of the Arzel\`{a}-Ascoli theorem are satisfied by bounding the spectral coefficients $\set{a^N_j}_{j=0}^N$ of $t \mapsto x^N$. A justification for this bound comes from Jackson's theorem,
%
\begin{theorem}[Jackson]\label{thm:Jackson}
Let $\dot y(t)$ be of bounded variation in $[-1, 1]$. Then, the coefficients,
$$a_j = \frac{1}{2} (2j+1)\int_{-1}^1 y(t)L_j(t) $$ \
satisfy the following inequality,
$$ \abs{a_j}  < \frac{6}{j^{3/2}\sqrt{\pi}}(M + V) $$
for all $j \ge 1$, where $L_j$ is the Legendre polynomial of degree $j$, $M$ is the least upper bound of $\abs{\dot y(t)}$ and $V$ is the total variation of $\dot y$ in $[-1, 1]$.
\end{theorem}
In other words, we select the Legendre PS method and further constrain the region $\X \times \U$ by a bound on the spectral coefficients of $\set{t \mapsto x^N}_{N=N_0}^\infty$ to ensure uniform equicontinuity. Details are provided in~\cite{RossGongBook, krg:convergence,kang:rate}.

It will be apparent in the next section that the practical implementation of these concepts is achieved by the spectral algorithm. Note that there are no coercivity type assumptions; hence, the theory does not require local uniqueness, and is applicable to problems with multiple optimal solutions.

Finally, consistency of the convergence requires that
$$ x^\infty(t) = x^*(t) \qquad  u^\infty(t) = u^*(t) \quad a.e. $$
where $t \mapsto (x^*, u^*) $ is an optimal solution of Problem $B$. This result is obtained by proving that $t \mapsto (x^\infty, u^\infty)$ is a feasible solution to Problem $B$ that achieves the optimal cost $E(x^*(-1), x^*(1))$. A key requirement for this proof is that $\lim_{N\to\infty}\pi^N = \pi^\infty$ is dense in $[-1, 1]$.  These results are encapsulated in the following theorem~\cite{RossGongBook, krg:convergence},
\begin{theorem}[Kang-Ross-Gong]\label{thm:kang}
Let
$\set{t \mapsto (x^{*N}, u^{*N})}_{N=N_0}^\infty$ be a sequence of optimal solutions generated by Problem $B^N$.  Then, there exists a subsequence, $\set{t \mapsto (x^{*n}, u^{*n})}_{n=n_0}^\infty$ and an optimal solution $t \mapsto (x^*, u^*)$ of Problem $B$ such that the following limits converge uniformly,
\begin{eqnarray}
\lim_{n \to \infty} x^{*n}(t)  &=& x^*(t)\nonumber\\
\lim_{n \to \infty} u^{*n}(t)  &=& u^*(t)\\
\lim_{n \to \infty} J^n[x^*(\pi^n), u^*(\pi^n)] &=& J[x^*(\cdot), u^*(\cdot)]\nonumber
\end{eqnarray}
\end{theorem}
The most difficult part of convergence theory is estimating the rate of convergence.  In~\cite{kang:rate}, Kang achieved this breakthrough for feedback linearizable systems. He showed that the rate of convergence of the cost function is $N^{1-2m_x/3}$; and, if the optimal control is $C^\infty$, then the convergence rate is faster than any given polynomial rate.

A number of additional results on convergence have been developed for the family of feedback linearizable systems with continuous or discontinuous control~\cite{kang:disc,kang:discontinuous-CDC}. One key strength of these convergence theorems is that they have been developed independent of the necessary conditions. Hence, when combined with the CMP, they provide a powerful tool for solving practical optimal control problems (see Problem $P$).

%

\subsection{Overview of Computational Implementation}
In any computational implementation, the best one can hope to achieve is to compute a solution to machine precision, $\epsilon_m > 0$ . This implies that we need to implement an algorithm to generate a solution up to some desired accuracy $\epsilon \ge \epsilon_m$.  From the results of the previous section, this implies that there exists an $N = N_\epsilon$ that solves the problem up to $\epsilon$ (including $\epsilon = \epsilon_m$) accuracy. In its computational implementation, the spectral algorithm~\cite{RossGongBook,knots,spec-alg,knots-automatic}, generates a sequence
$$\set{t \mapsto (x^{*N}, u^{*N})}_{N=N_0}^{N_\epsilon}$$
such that,
$$\norm{x^{*N_\epsilon}(\cdot) - x^*(\cdot)} \le \epsilon \quad and \quad  \norm{u^{*N_\epsilon}(\cdot) - u^*(\cdot)} \le \epsilon $$
where appropriate Sobolev norms, $W^{m_x,p}$ and $W^{m_u, p}$, are used for measuring errors in the state and control trajectories.  The conditions for Theorem \ref{thm:kang} are met by generating the spectral coefficients for $x^{*N}(t)$ using the relation,
\begin{eqnarray*}
a_j^N =(j + 0.5) \sum_{i=0}^N L_j(t_i) w_i x^*_i \qquad j = 0, \ldots, N
\end{eqnarray*}
and imposing a Jackson stopping criterion,
$$\abs{a^{N_\epsilon}_{N_\epsilon}} \le \delta_x$$
where $\delta_x > 0$ is some given tolerance and $x^*_i, \ i = 0, \ldots, N $ is the solution to Problem $B^N$ obtained by selecting $\pi^N$ as a Gauss-Lobatto grid for finite horizon problems and a Gauss-Radau grid for infinite-horizon problems (see Remark \ref{rm:summary}). This concept is part of a broader spectral algorithm~\cite{RossGongBook,spec-alg,guess-free,DIDO_2007} wherein additional conditions are imposed such as a dual stopping condition for the adjoint covector derived by an application of the CMP. 
When all convergence criteria are met, the algorithm terminates.

\begin{figure}[h!]
      \centering\scalebox{0.5}
      {\includegraphics[scale=0.65]{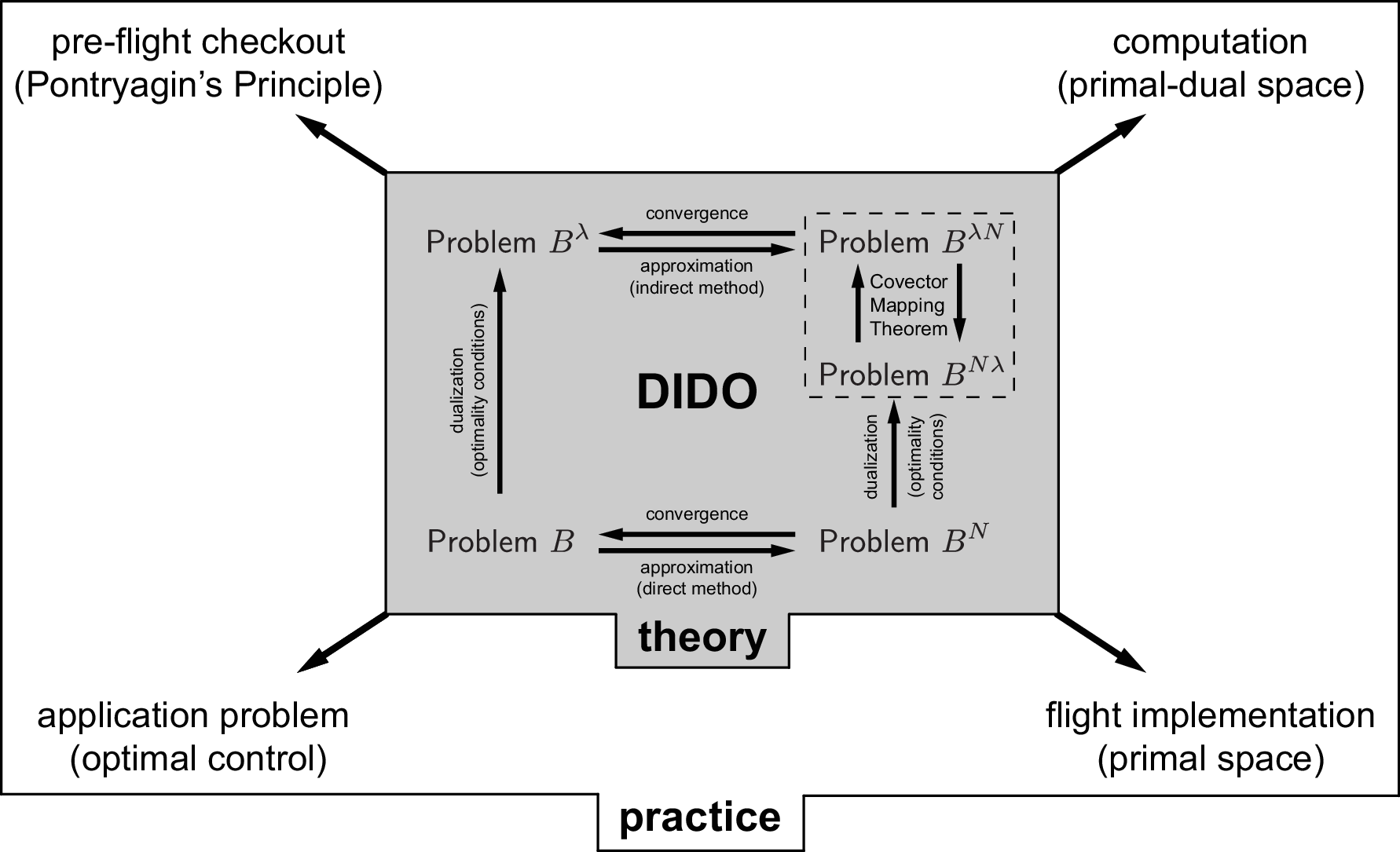}}
      \caption{\textsf{Illustration of the close connection between software, theory and practice of PS optimal control.}}
      \label{fig:CMP-layered}
   \end{figure}

Currently, there are two major implementations of the spectral algorithm. One is in OTIS, a NASA software package written in FORTRAN with libraries to support aerospace trajectory optimization problems~\cite{otis,paris:PS}.  The other implementation is in DIDO, a MATLAB based general purpose commercial optimal control software package~\cite{DIDO_2007}. DIDO also holds the distinction of being the first implementation of PS optimal control. It contains many of the advancements noted in this paper including an implementation of pseudospectral knots~\cite{knots}.  PS knots enable the practical implementation of discontinuous controls as well as jumps in the state variables, i.e. hybrid-type optimal control problems~\cite{acc:hybrid,hybrid:jgcd,stevens:cycler,josselyn:aerocapture}.  All of the ground and flight implementations of PS control have been performed using DIDO. Fig.~\ref{fig:CMP-layered}
%
%
shows how the software, theory, and flight implementation blend together.  It is clear that theory, computation and practice are intimately linked.

\section{Ground and Flight Implementations}
\begin{figure}[h!]\centering\scalebox{0.55}
      {\includegraphics[scale = 0.6]{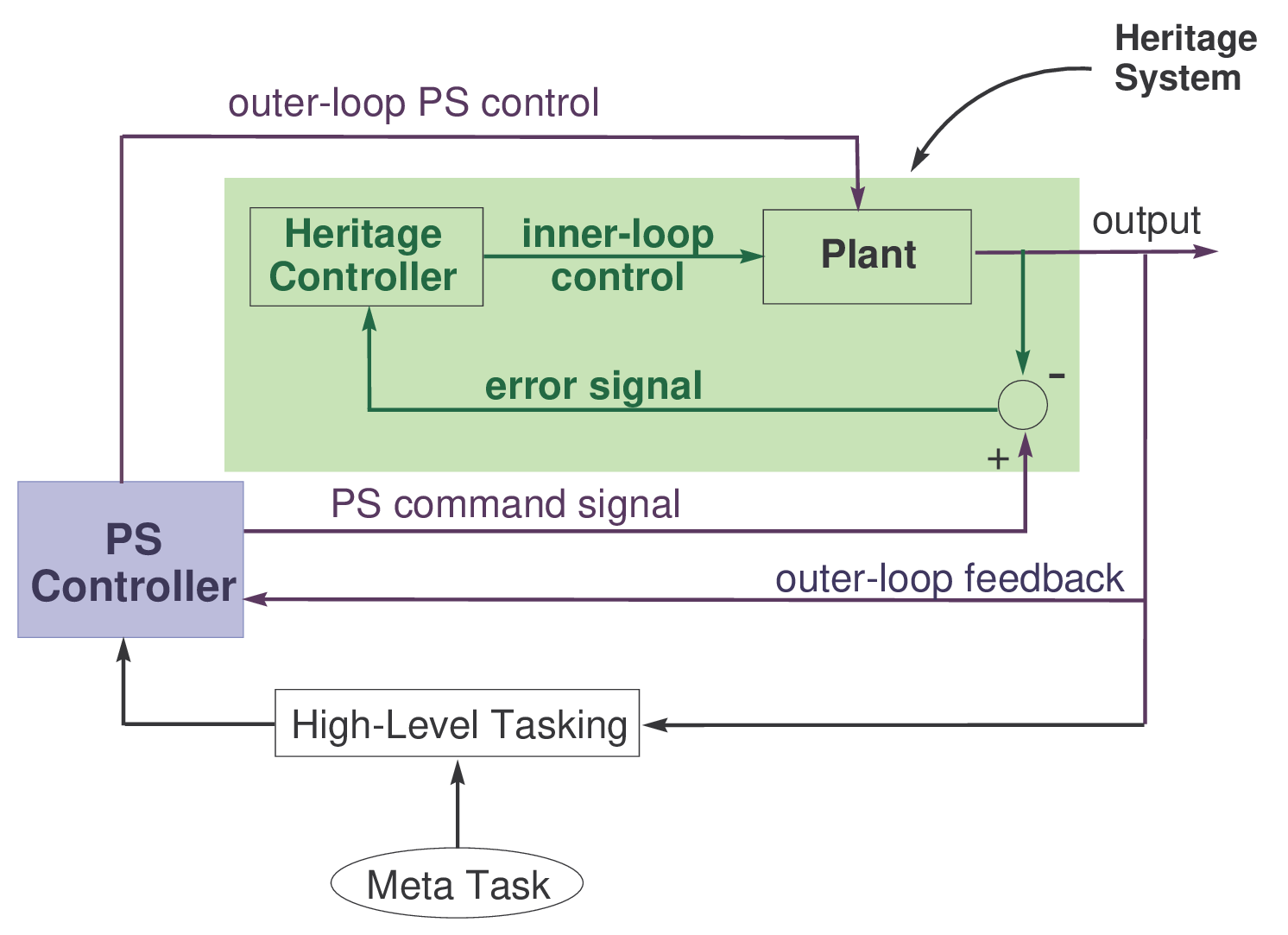}}
      \caption{\textsf{Implementation of PS controls as ``strap-on'' enhancements to heritage systems.}}
      \label{fig:PS-implement-block-diag}
   \end{figure}

In 2007, a special session at the IEEE Conference on Decision and Control highlighted the then emerging role of PS control in military and industrial applications~\cite{CDC:session}. Early applications of PS concepts were limited to open-loop simulations~\cite{rea,stanton}.  A series of rapid breakthroughs~\cite{unifed-CDC,flat,issues,unified-wObserver} with ground and flight implementations reshaped the role of PS control towards a design and implementation of outer loops as illustrated in Fig.~\ref{fig:PS-implement-block-diag}.
%
%
PS feedback control is obtained via real-time computation of optimal control~\cite{RossGongBook,acc:stability,PS:Radau,williams:RHC-TN,jgcd:c-pi}.
In this section we describe signature ground and flight implementations that illustrate some key features of industrial importance.

\subsection{Ground Implementation Example}
\begin{figure}[h!]
      \centering\scalebox{1}
      {\includegraphics[scale=1]{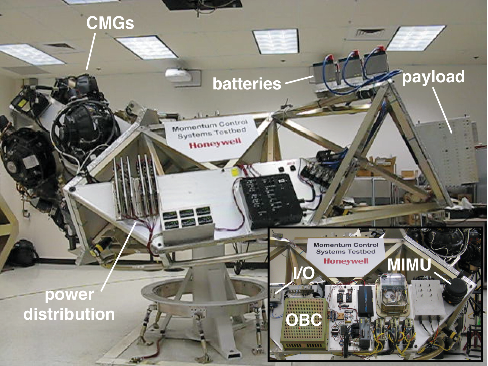}}
      \caption{\textsf{Honeywell's industrial-grade momentum control system.}}
      \label{fig:MCS_testbed}
   \end{figure}

While many ground implementations of PS control have been performed~\cite{RossGongBook}, a particular execution carried out at Honeywell illustrates its industrial use.  A 3000 lbs industrial-grade momentum control system (MCS) assembly (shown in  Fig.~\ref{fig:MCS_testbed})
%
%
floats on a powerful airbearing in a large test facility of approximately 1100 square feet. The testbed contains six control moment gyroscopes (CMGs) mounted in a pyramid configuration. Because the CMGs are mounted off-center, the system is balanced to align the testbed center of gravity with the center of rotation. This is accomplished in two main stages~\cite{honeywell}: (\emph{i}) manual balancing and (\emph{ii}) automatic balancing.  Manual balancing is carried out by positioning of the lead acid batteries (which provide power for the testbed) and stackable payload weights on a grid with 25-mm spacing. 
The second mass balancing stage utilizes a set of three prismatic actuators with very fine position control to locate a series of masses so as to precisely drive the center of mass to the air bearing pivot.

\begin{figure}[h!]
      \centering\scalebox{1}
      {\includegraphics[scale=0.5]{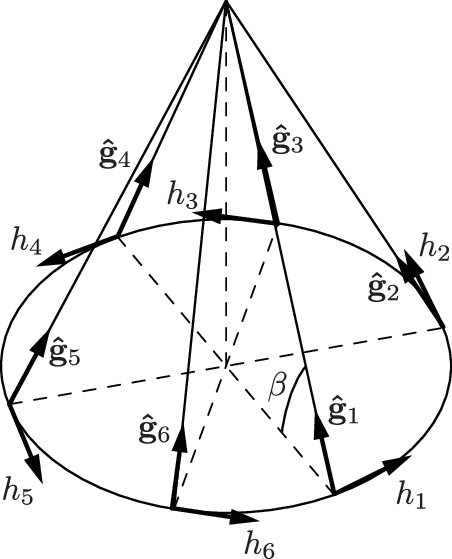}}
      \caption{\textsf{Schematic of six CMGs arranged as part of a skew array.}}
      \label{fig:skew_array}
   \end{figure}

The CMG units are arranged in a pyramidal fashion about a base circle (see Fig.~\ref{fig:skew_array}) and produce momentum along  $\bh_i$.
%
%
Torque is produced by a rotation of the gimbals around $\mathbf{\hat g}_i$ according to the following equation,
$$\tau = \mathbf{A}(\delta)\dot\delta  $$
where $\dot\delta$  is the gimbal rotation rate.

Since the CMG system is inherently over actuated, the torque command vector must be converted to a gimbal rate command,  $\dot\delta_c$, to drive the CMGs
The simplest steering logic is to apply the Moore-Penrose pseudoinverse to generate,
$$ \dot\delta_c = \mathbf{A}^T(\delta) [\mathbf{A}(\delta) \mathbf{A}^T(\delta)]^{-1}\tau  $$
A challenge problem related this type of control allocation is the so-called singularity problem. This condition occurs when the matrix product, $\mathbf{A}(\delta) \mathbf{A}^T(\delta)$ , becomes degenerate and the pseudoinverse steering can no longer be computed. In the singular condition, the geometric configuration of the CMGs is such that the net CMG torque is available in all but one direction. The direction in which a CMG torque cannot be produced is called the singular direction~\cite{Wie_2004}. On a singular surface (a collection of singular directions), it is not possible to change the system momentum vector regardless of the control input. In order to control the system, therefore, it is necessary to design maneuvers so that the momentum trajectory avoids the singular surface.

The singularity problem is addressed in heritage industrial applications by building conservatism into the momentum envelope. This approach specifically avoids singular regions by reducing the momentum envelope so that the singularities can be excluded from the momentum space used to operate the system. The conservatism in the reduced momentum envelope generates an increase in size, weight and power of the system to meet a given rate requirement. Singularity avoidance techniques and singularity escape algorithms can also be used to generate steering laws that do not become degenerate when a singularity surface is approached. One simple example of such a singularity-robust steering logic is given by~\cite{Bedrossian_1990}:
$$ \dot\delta_c = \mathbf{A}^T(\delta)[\mathbf{A}(\delta)\mathbf{A}^T(\delta) + \lambda \mathbf{I}]^{-1} \tau$$
where $\mathbf{I}$ is the identity matrix and $\lambda$ is a small positive constant. The addition of $\lambda \mathbf{I} $  maintains the linear independence of the rows of matrix $\mathbf{A}$ so that $\det(\mathbf{A}(\delta)\mathbf{A}^T(\delta)) \ne 0 $. When the value of $\det(\mathbf{A}(\delta)\mathbf{A}^T(\delta))$  is large, the singularity-robust steering logic behaves nearly identically to the Moore-Penrose pseudoinverse. As $\det(\mathbf{A}(\delta)\mathbf{A}^T(\delta)) \to 0$, the steering logic prevents the occurrence of singular configurations by perturbing the momentum vector around the singular states.

\begin{figure}[h!]
      \centering 
      {\includegraphics[scale=0.5, clip]{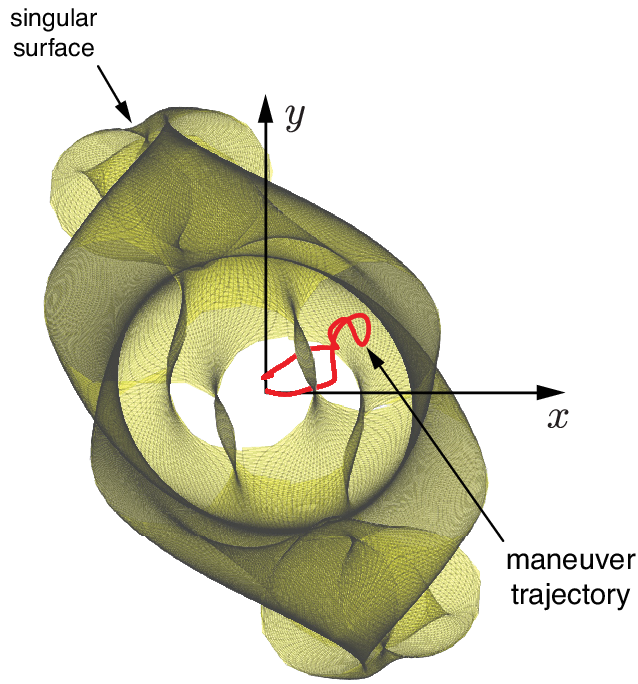}}
      \caption{\textsf{Singularity avoidance maneuvering in momentum space.}}
      \label{fig:momentum_plot}
   \end{figure}

A radically different approach to overcoming the singularity challenge problem is to design and implement optimal solutions.  The the feasible region is defined in a manner that excludes the singularity region. For instance, for a typical skew-type array, the mixed-state-control path constraint,
%
\begin{equation*}
\left[h_x -
\left(
\begin{split}
 &\frac{\cos\beta(u_x\cos\beta - u_z\sin\beta)}{\sqrt{1 - (u_x\sin\beta + u_z\cos\beta)^2}} +
\frac{u_x}{\sqrt{1 - (u_y\sin\beta + u_z\cos\beta)^2}} +\\
&\qquad \qquad \frac{\cos\beta(u_x\cos\beta + u_z\sin\beta)}{\sqrt{1 - (u_z\cos\beta - u_x\sin\beta)^2}} +
\frac{u_x}{\sqrt{1 - ( u_z\cos\beta - u_y\sin\beta)^2}}
\end{split}
\right)
\right]^2 \ge \alpha
\end{equation*}
avoids an $\alpha > 0 $ region specified by the designer.
This concept is similar to the notion of obstacle avoidance in robotics~\cite{infotech,cascio_2009,hurni_2010,hurni-main}.   Thus, the optimal control framework is used to design and implement feasible solutions with optimality as a side benefit. The optimality criterion can be anything sensible such as electric power or other practical measures of performance.  A practical PS realization of these concepts for the Honeywell testbed is shown in Fig.~\ref{fig:momentum_plot}.  The PS concept was implemented as part of a strap-on implementation (see Fig.~\ref{fig:PS-implement-block-diag}) to Honeywell's 6-CMG skew array with DIDO serving as the PS controller.  It is apparent from Fig.~\ref{fig:momentum_plot} that the maneuver avoids the singular surface regions.
%
%

%
%

\subsection{Flight Implementation Example}
\begin{figure}[h!]
      \centering\scalebox{1}
      {\includegraphics[scale=0.5]{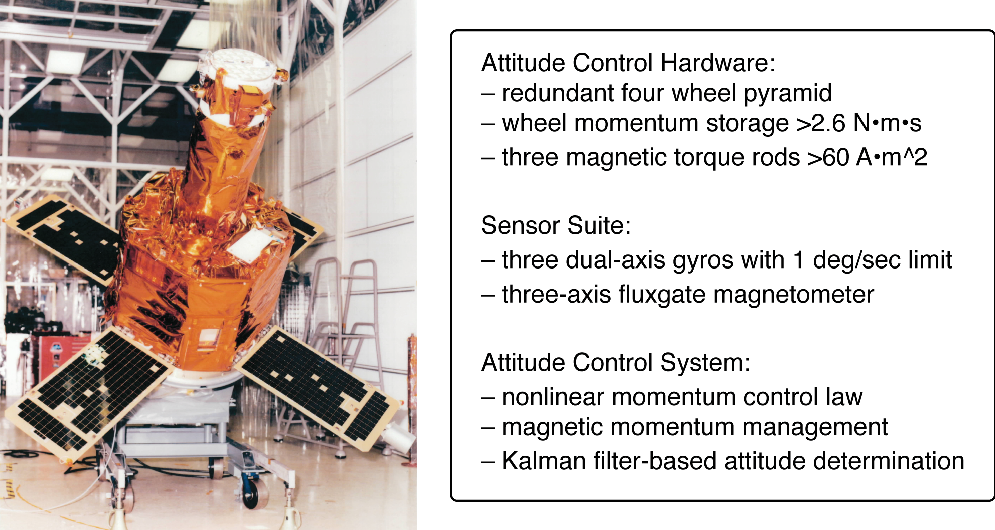}}
      \caption{\textsf{The TRACE spacecraft undergoing a pre-launch checkout at NASA.}}
      \label{fig:TRACE}
   \end{figure}

\begin{figure}[h!]
      \centering\scalebox{1}
      {\includegraphics[scale=0.5]{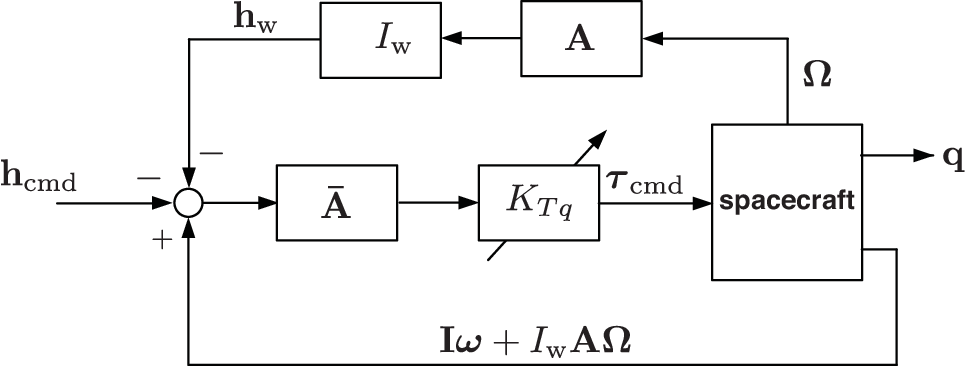}}
      \caption{\textsf{Block diagram of TRACE's heritage closed-loop momentum control system.}}
      \label{fig:TRACE_momentum_ACS}
   \end{figure}

As noted in Section 1, the first flight implementation of PS control was in 2006. Since then there have been a number of other flight implementations and more are on a path for operational use~\cite{RossGongBook}. The 2010 flight implementation of PS optimal control performed on NASA's space telescope, TRACE (see Fig.~\ref{fig:TRACE}), is discussed in this section to illustrate the close connection between mathematical theory and flight operations.
%
%
TRACE is responsible for measuring fine scale magnetic features of the solar surface and corona. To maneuver between points of interest, TRACE is outfitted with a reaction wheel-based attitude control system.
Three magnetic torque coils are also employed for momentum management.
Momentum transfer is accomplished by using a closed-loop momentum control system (see Fig.~\ref{fig:TRACE_momentum_ACS})
%
%
to control the angular rates of the individual reaction wheels.

\begin{figure}[h!]
      \centering\scalebox{1}
      {\includegraphics[scale=0.6]{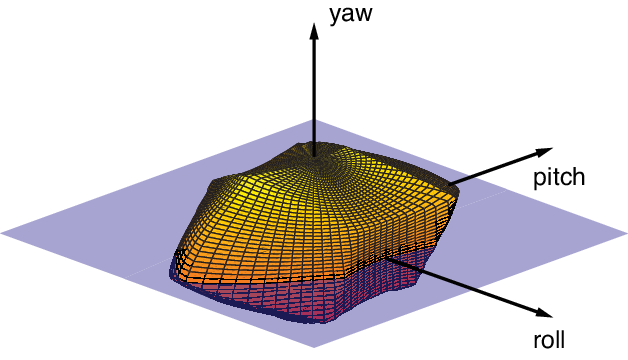}}
      {\includegraphics[scale=0.6]{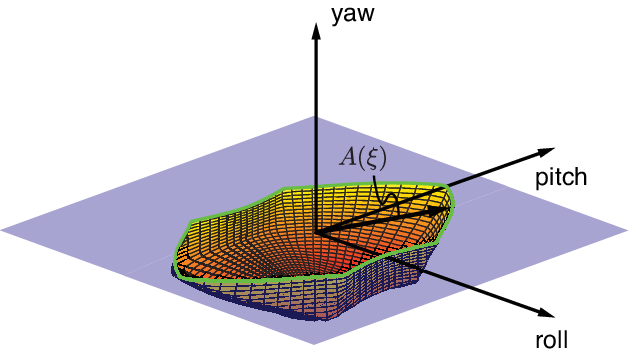}}
      \caption{\textsf{Agilitoid for the TRACE spacecraft (left) and its cutaway view (right).}}
      \label{fig:TRACE_envelope}
   \end{figure}

The agility capability of a spacecraft can be visualized in terms of an ``agilitoid;'' see Fig.~\ref{fig:TRACE_envelope}.
%
%
The agilitoid is generated by mapping the momentum-to-inertia ratio,
$$ A(\xi) = \frac{h(\xi)}{I_\xi}  $$
over a $2\pi$ steradian, where $h(\xi)$ is the angular momentum of the spacecraft about an arbitrary axis, $\xi$, and $I_\xi$ is the moment of inertia of the spacecraft along $\xi$. Fig.~\ref{fig:TRACE_envelope}
clearly shows that the agility of TRACE is greatest around the pitch axis and smallest around the roll and yaw axes. Thus, a minimum-time space maneuver should be able to exploit this nonspherical volume to generate new, counter-intuitive maneuvers similar to the Brachistochrone problem. This problem can be formulated as the shortest-time maneuver (STM),
\begin{equation}
\textsf{(STM) }\left\{   \begin{array}{lrll}
                \textsf{Minimize }  & J [\bx(\cdot), \mathbf{\tau}_{cmd}, t_f] &=t_f \\
                \textsf{Subject\ to}
& \dot \bx(t) &= \left \{ \begin{array}{l}
        \frac{1}{2}\mathbf{Q}(\mathbf{\omega})\mathbf{q} \\
        \mathbf{\Gamma}^{-1} \left [
                \begin{array}{c}
                                        - \mathbf{\omega} \times \left( \mathbf{J} \mathbf{\omega} +
                        \sum_{i=1}^4 \mathbf{a}_i I_{{ w},i} \Omega_{{ w},i} + \mathbf{a}_i I_{{ w},i} \mathbf{a}_i^T \mathbf{\omega} \right) \\
                                        \mathbf{C}\bx_{{filt}} + \mathbf{D}\mathbf{\tau_{{cmd}}}
                \end{array}
        \right ] \\
        \mathbf{A}\bx_{{filt}} + \mathbf{B}\mathbf{\tau_{{cmd}}} \\
\end{array} \right \} \\
& \bx(t_0) &=  \left [\mathbf{e}_0\sin(\frac{\phi_0}{2}), \cos(\frac{\phi_0}{2}),\mathbf{\omega}_0, \Omega_0, \bx_{{filt},0})\right]^T\\
& \bx(t_f) &=  \left [\mathbf{e}_f\sin(\frac{\phi_f}{2}), \cos(\frac{\phi_f}{2}),\mathbf{\omega}_f, \Omega_f, \bx_{{filt},f}) \right]^T \\
&     \vert \vert \mathbf{q}(t) \vert \vert &= 1 \\
&       \vert \omega_i(t) \vert  &\leq \omega_{\max} ,\quad i=1,\ldots,3 \\
&     \vert \tau_{{cmd},i}(t) \vert &\leq \tau_{{cmd,max}},\quad i=1,\ldots,4 \\
&     \vert I_{{w},i}\Omega_{{w},i}(t) \vert &\leq I_{{w},i} \Omega_{\max},\quad i=1,\ldots,4 \\
&     \tau^L & \leq I_{{w},i}\dot \Omega_{{w},i}(t)  \leq \tau^U,\quad i=1,\ldots,4\\
        \end{array} \right.
        \label{eq_13}
\end{equation}
The spacecraft dynamics include the quaternion kinematics, the multi-body dynamics, as well as the dynamics of a bank of reaction wheel command-shaping filters. The command-shaping filters are used to suppress the excitation of spacecraft structural modes.  In order to ensure the implementability of the optimal solution on flight hardware, the problem formulation also includes a number of practical constraints on the commandable reaction wheel torques, rate constraints to avoid saturating the onboard rate gyros, electrical power constraints and momentum limits.  All these engineering and operational constraints were implemented as nonlinear path constraints of the form, $\bh^L \le \bh(\bx(t), \bu(t)) \le \bh^U$; see Problem $P$ discussed in Section 2.

\begin{figure}[h!]
      \centering\scalebox{1}
      {\includegraphics[scale=0.4]{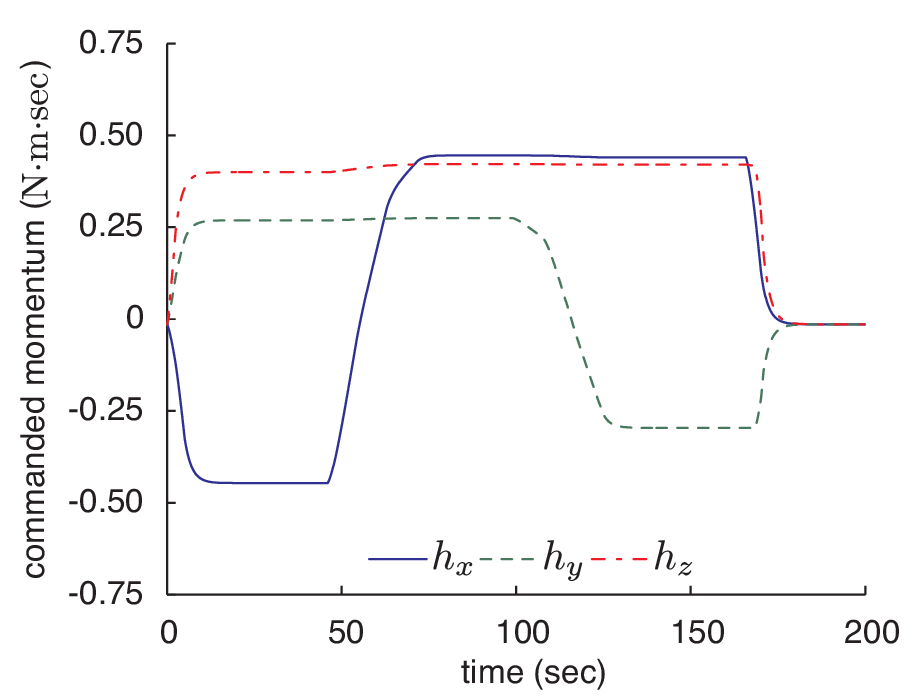}}
      {\includegraphics[scale=0.4]{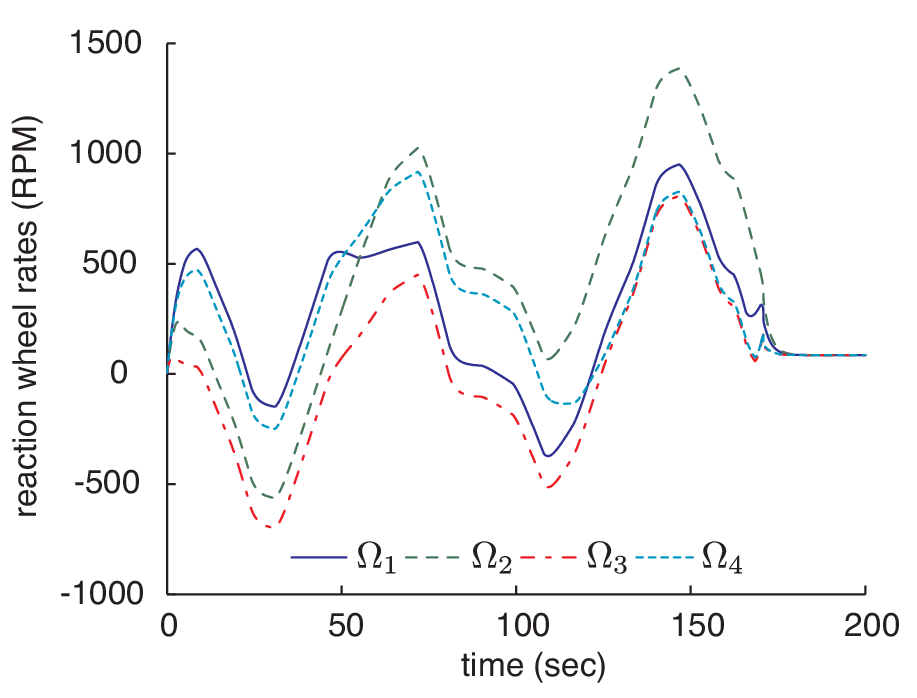}}
      \caption{\textsf{Initial pre-flight STM solution for TRACE.}}
      \label{fig:MODE1}
   \end{figure}

\begin{figure}[h!]
      \centering\scalebox{1}
      {\includegraphics[scale=0.4]{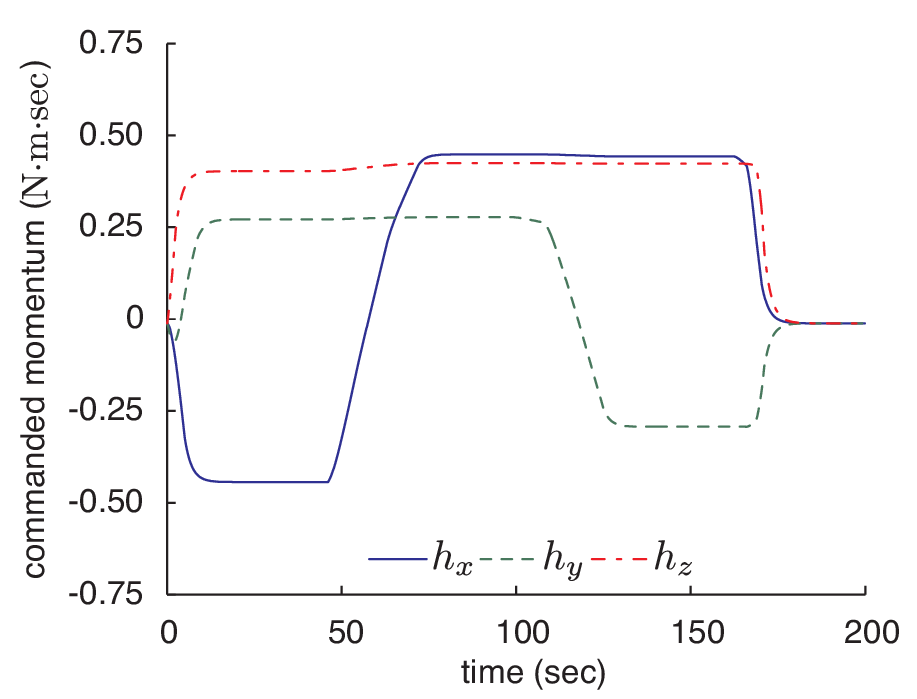}}
      {\includegraphics[scale=0.4]{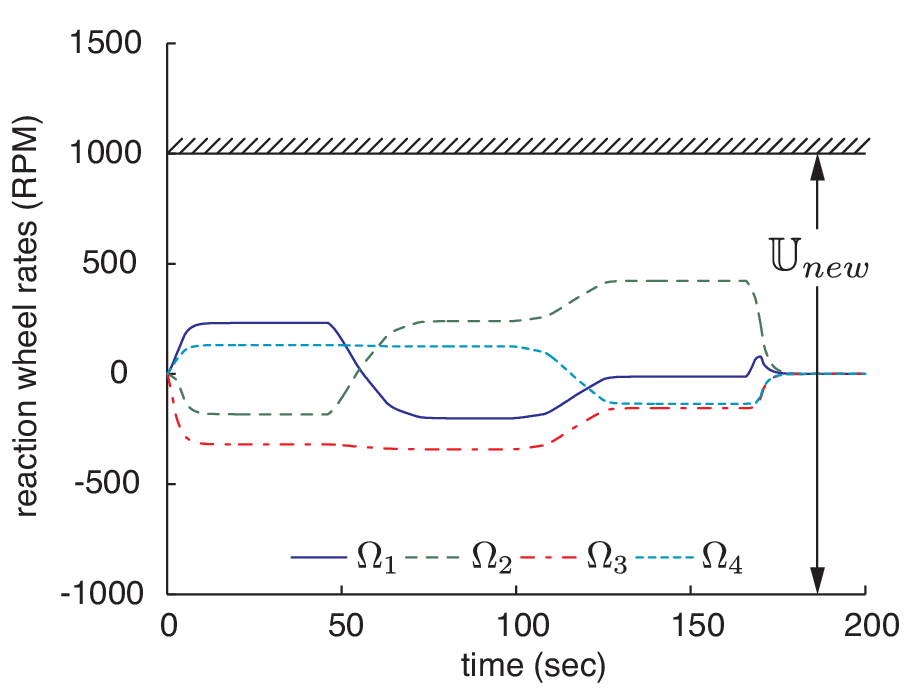}}
      \caption{\textsf{Proposed STM implementation for TRACE after pre-flight checkout (as illustrated in Fig.~\ref{fig:CMP-layered}).}}
      \label{fig:MODE2}
   \end{figure}

\begin{figure}[h!]
      \centering\scalebox{1}
      {\includegraphics[scale=0.4]{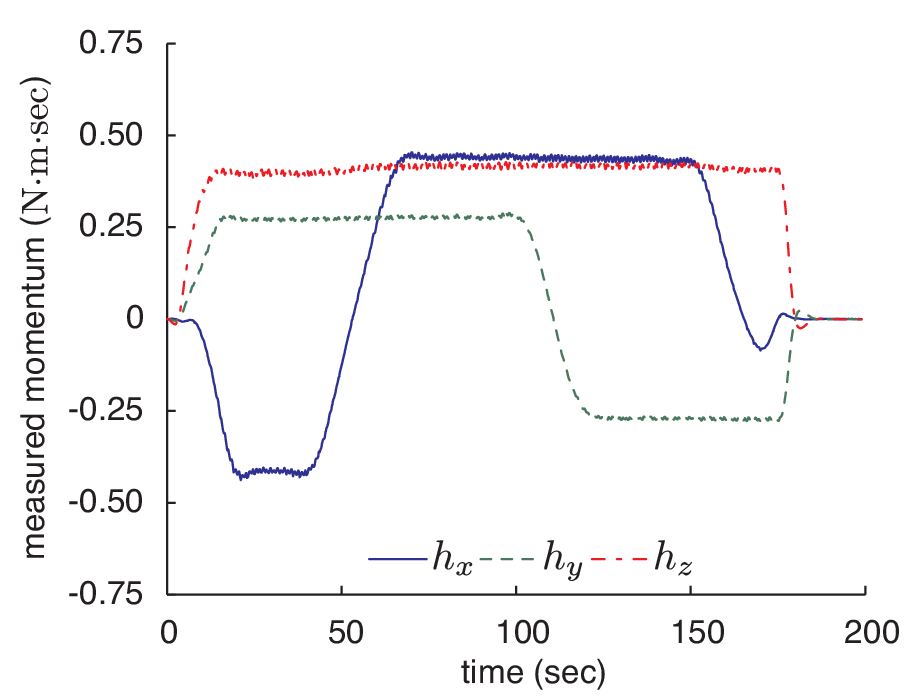}}
      {\includegraphics[scale=0.4]{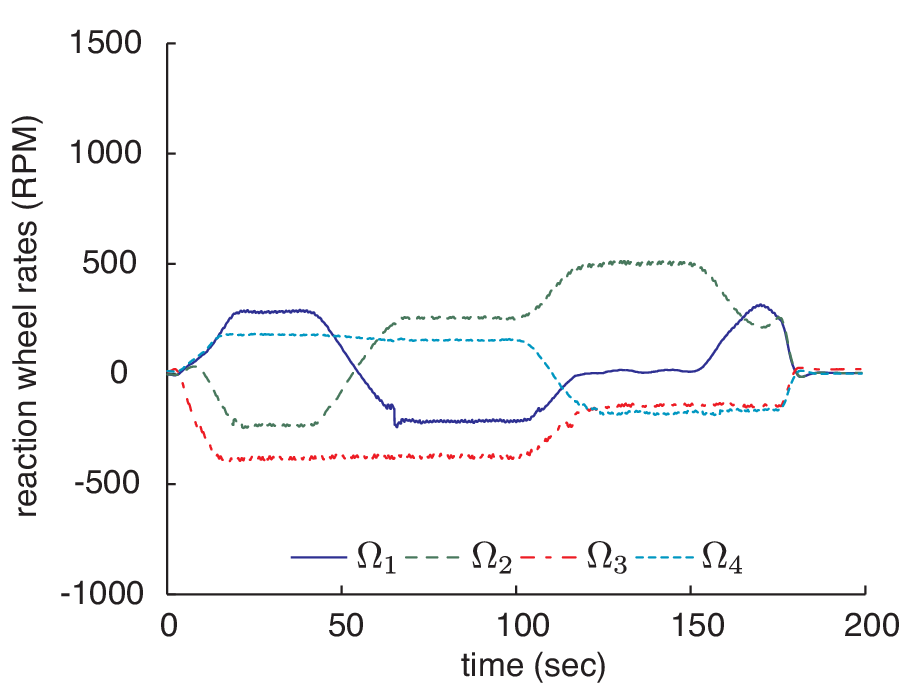}}
      \caption{\textsf{TRACE flight results for the implementation proposal of Fig.~\ref{fig:MODE2}.}}
      \label{fig:flight}
   \end{figure}

In order to interface with TRACE's heritage control system, momentum command trajectories were generated using the relation,
$$ \bh_{cmd} = I_w \mathbf{A} \int \tau_{cmd}\ dt $$
where $\mathbf{A}$ denotes the transformation matrix that maps the reaction wheel momenta from the individual actuator frames to the spacecraft body-fixed frame. Fig.~\ref{fig:MODE1}
%
%
shows the candidate optimal solution obtained via DIDO.  As part of the pre-flight checkout activities (see Fig.~\ref{fig:CMP-layered}) engineers felt that the reaction wheel rates, while meeting all the specifications, were a little unusual in the sense that they all changed speeds in unison with one another.  They speculated a ``better'' solution existed.  In order to explore this conjecture, the acceptable control space was reduced from, $\U_{old} = \set{\abs{\Omega_i} \le 2500\ RPM, \ i = 1, 2, 3, 4 } $ to $\U_{new} = \set{\abs{\Omega_i} \le 1000\ RPM, \ i = 1, 2, 3, 4 } $ as a means to select a different subsequence for convergence (see Theorems \ref{thm:kang} and \ref{thm:AA}). The resulting solution is shown in
Fig.~\ref{fig:MODE2}.
%
%
Note that the box constraints are still inactive and the momentum solution remains the same. This implies that the spectral algorithm was able to choose a different subsequence for convergence in accordance with Theorem \ref{thm:kang}.  A flight implementation of this result is shown in Fig.~\ref{fig:flight}.
%
%
It is clear that the telemetry data are in close agreement with the PS predicted results indicating the validity of the entire process.

\section{Embedded Optimal Control}
\begin{figure}[h!]
      \centering\scalebox{1}
      {\includegraphics[scale=0.5]{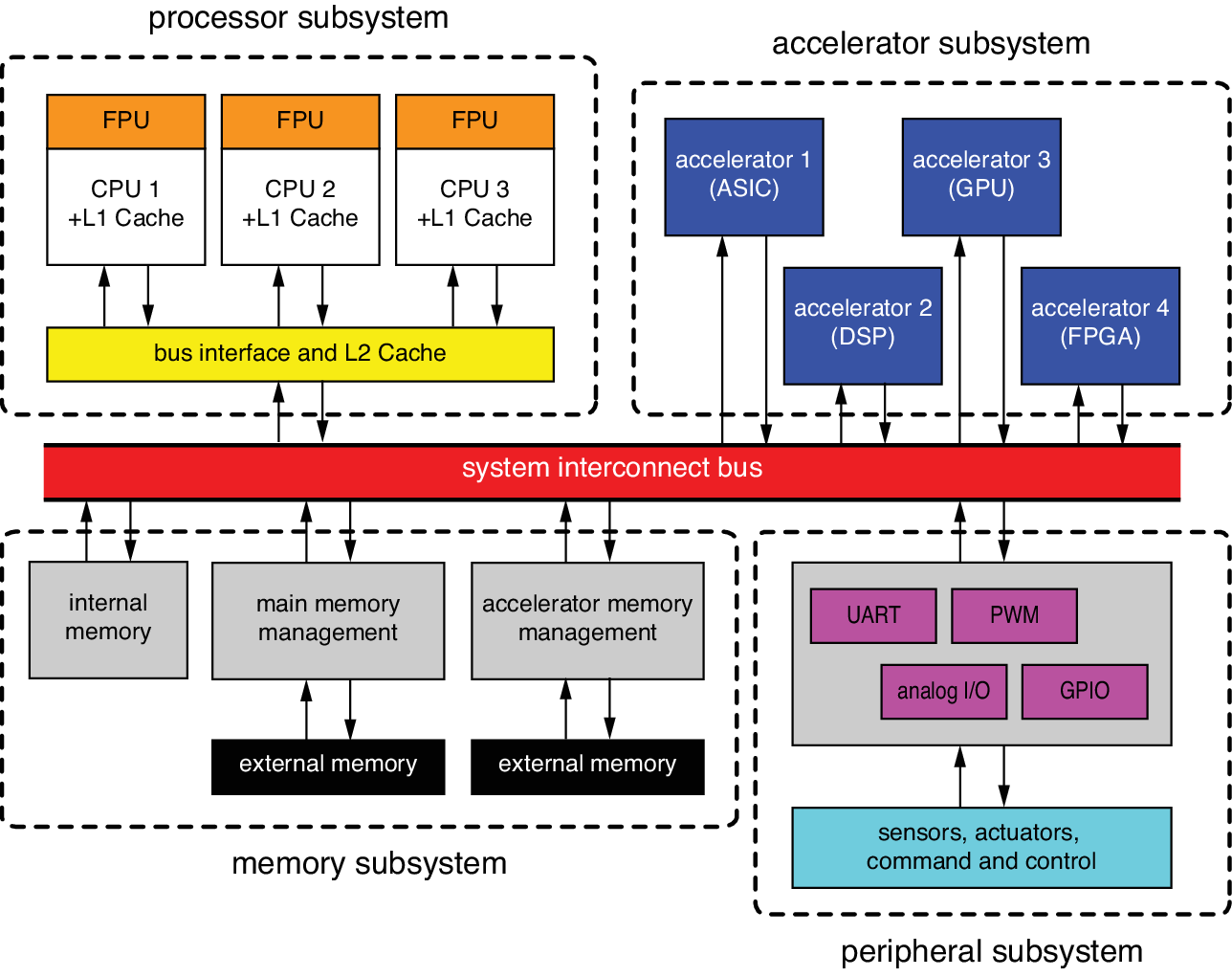}}
      \caption{\textsf{Architecture for an embedded PS optimal control system.}}
      \label{fig:EMBEDDED_arch}
   \end{figure}

The spectral algorithm can be embedded in an assemblage of special purpose digital circuits for the principal means of performing efficient and fast computations of PS controls. In sharp contrast to a hardware implementation of an industry standard PID controller that may only need one small field programmable gate array (FPGA), a PS controller requires a more sophisticated architecture as shown in Fig.~\ref{fig:EMBEDDED_arch}.
%
%
 This is because the goal of the embedded spectral algorithm is to solve Problem $B$ to a given $\epsilon > 0 $ tolerance within a given time budget, $\tau_b$.  It can be shown\cite{acc:stability,jgcd:c-pi} that the quantity $\tau_b$ varies inversely as the Lipschitz constant of the vector field, ${\bf f}(\bx,\bu,t)$, with respect to $\bx$.  Hence, the architecture of Fig.~\ref{fig:EMBEDDED_arch} can be appropriately tailored for specific plants (aircraft, spacecraft, robotic systems, etc). Regardless, for any given system, the spectral algorithm must perform well-defined tasks such as the computation of the grid points, $t_j, \ j = 0, \ldots, N$ that define $\pi^N$, the differentiation matrix, $D_{ik}$, integration weights, $w_i$, etc.  These tasks can be properly parsed within the architecture shown in Fig.~\ref{fig:EMBEDDED_arch} to produce a generic embedded PS control solution for Problem $B$ (and hence Problem $P$).

\subsection{Hardware Architecture}
An embedded PS optimal control system comprises of four main subsystems: processor, accelerator, memory and a peripheral subsystem that are connected by a core system interconnect bus.
The main processor subsystem contains one or more central processing unit~(CPU). Each CPU core also contains an integrated floating point unit (FPU). A multi-core processor yields distinct advantages since it is possible for some tasks to be distributed across more than one core to improve computational throughput. Thus, the challenge is to find ways to implement the spectral algorithm that best leverages the inherent capabilities of the various computational subsystems.  One aspect of distributing these computations is through the accelerator subsystem. Similar to an FPU, hardware accelerators help reduce the computational time by using dedicated hardware to perform repetitive tasks more quickly than is possible by a software instantiation running on a main processing core~\cite{Moon_2005}. One example of such a hardware accelerator is the Digital Signal Processor (DSP) chip which is extremely useful in the computation of the differentiation and integration matrices (see equations \ref{eq:D-def-Wt} and \ref{eq:weights}).

The inclusion of Graphical Processor Units (GPUs) as part of the hardware accelerator subsystem allows some of the necessary operations for PS optimal control computation to be performed more rapidly. In addition, sparse matrix multiplication and LU decomposition have been designed for low-power embedded FPGA hardware that can exceed the performance of optimized routines running on much larger desktop computers~\cite{fpga}. In fact, due to their user programmable nature, FPGAs are an essential tool in the development and testing of embedded spectral algorithms. Implementing and testing various aspects of the spectral algorithms on FPGA accelerators allows the ideal embedded hardware architecture to be iterated and ultimately configured.

Besides outfitting the embedded system with sufficient memory resources for solution computation, the memory subsystem must also manage data transfers between the main processors and the various hardware accelerators. Thus, the architecture of the interconnect bus (see Fig.~\ref{fig:EMBEDDED_arch}) can have a great influence on the overall performance of the embedded PS controller. Adopting System-on-Chip (SoC) or Computer-on-Module (CoM) paradigms, in which multiple components of the embedded system architecture are integrated onto a single chip, helps manage this problem.

To communicate with the outside world, the embedded optimal control system is also outfitted with a peripheral interface system. The peripheral interface subsystem implements a variety of bus architectures for communicating with hardware external to the embedded PS control hardware. These architectures can include serial protocols such as RS-232/485, I$^2$C, SPI, USB, CAN, or Ethernet. Together with various analog and digital I/O as well as pulse-width-modulation (PWM) drivers, the system connects with sensors, actuators and external command and control systems as illustrated in Fig.~\ref{fig:PS-implement-block-diag}.

\subsection{Streamlining Embedded System Performance}

\begin{figure}[h!]
      \centering\scalebox{1}
      {\includegraphics[scale=0.4]{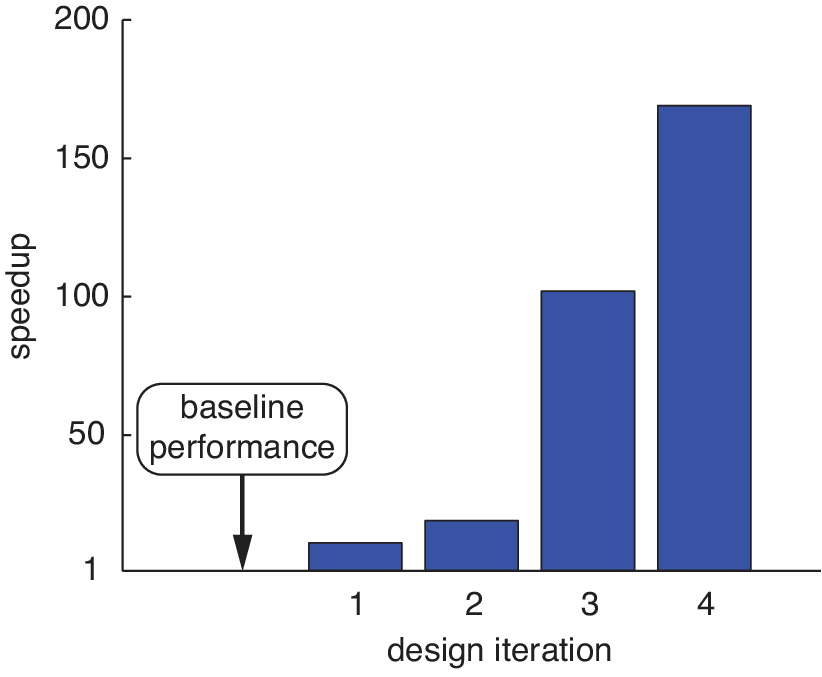}}
      \caption{\textsf{Computational speedups achieved by various parsings of the spectral algorithm for the shortest-time space maneuver.}}
      \label{fig:EMBEDDED_speedup}
   \end{figure}

Recent advances in the development and implementation of embedded PS controllers are illustrated in Fig.~\ref{fig:EMBEDDED_speedup}.
%
%
Fig.~\ref{fig:EMBEDDED_speedup} shows the typical computational speedup that can be achieved as a result of iterating on the design of the embedded spectral algorithm as well as refining the distribution of the algorithm among the various computational units shown in Fig.~\ref{fig:EMBEDDED_arch}. The first design iteration illustrates the performance improvement obtained by transforming the baseline software into a directly executable machine code. Subsequent design iterations involve refinements made to the baseline that leverage the specific capabilities of the subsystem units (see Fig.~\ref{fig:EMBEDDED_arch}). Computational speedups of over 175 times have been obtained without even exhausting all of the capabilities of the architecture illustrated in Fig.~\ref{fig:EMBEDDED_arch}.  Thus, additional modifications are possible to further improve the computational performance of the spectral algorithm.

\begin{figure}[h!]
      \centering\scalebox{1}
      {\includegraphics[scale=0.4]{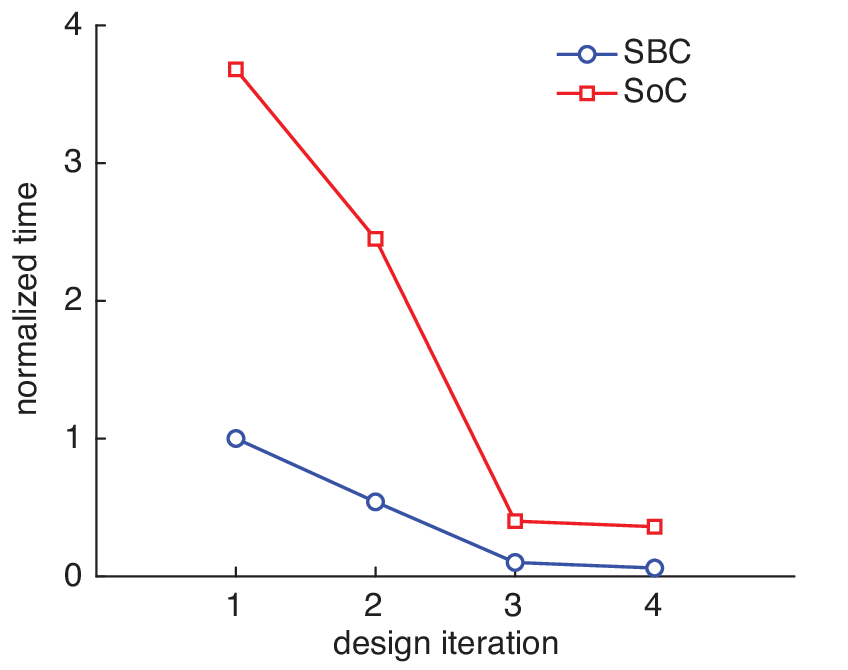}}
      \caption{\textsf{Normalized solve times for SBC- and SoC-type architectures.}}
      \label{fig:EMBEDDED_solve_time}
   \end{figure}

\begin{figure}[h!]
      \centering\scalebox{1}
      {\includegraphics[scale=0.4]{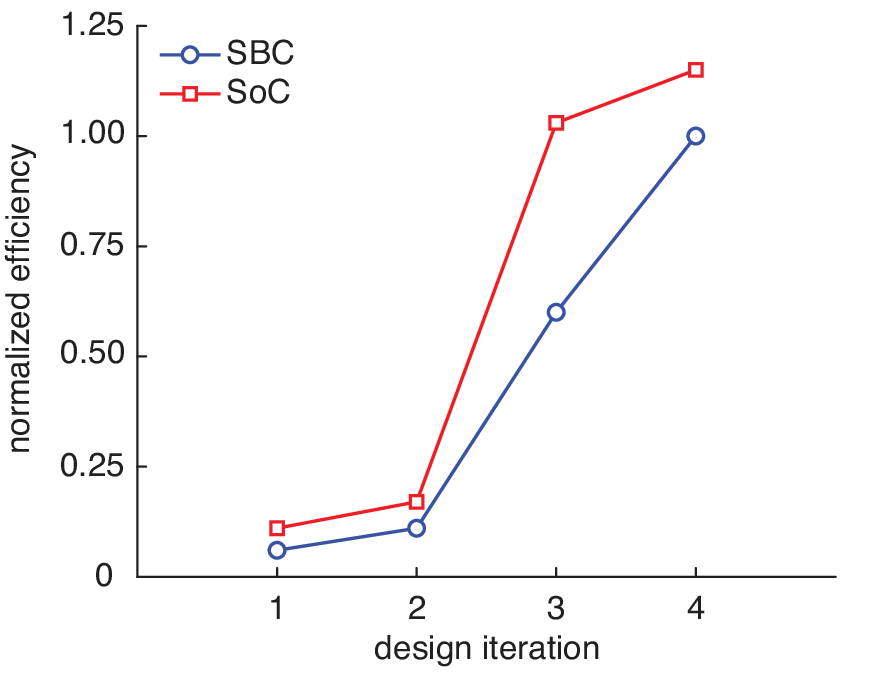}}
      \caption{\textsf{Normalized performance per watt for SBC- and SoC-type architectures}}
      \label{fig:EMBEDDED_perf_per_watt}
   \end{figure}

To highlight the influence of the architecture design on the overall performance of the embedded optimal control system, the relative performance of two different prototype architectures is compared in Fig.~\ref{fig:EMBEDDED_solve_time}.
%
%
%
It is clear that the relative computational performance different architectures varies significantly. A Single Board Computer (SBC)-like architecture solves Problem $B$ approximately three times faster than a SoC-type architecture in  design iteration 1. This is primarily due to differences in the main processing and hardware acceleration subsystems. The overall performance of each architecture can be significantly improved by identifying and isolating architecture-specific computational bottlenecks and re-addressing the algorithm parsing problem in an iterative fashion. As evident in Fig.~\ref{fig:EMBEDDED_solve_time}, this process can quickly diminish the absolute advantage of a given architecture in terms of the metric of solve time. Thus, a second metric for comparing the performance of various architectures is the performance per watt. This metric describes the rate of computation per watt of power consumed and can be an important deciding factor in certain applications (e.g. space) where power consumption is always at a premium. This analysis shows (see Fig.~\ref{fig:EMBEDDED_perf_per_watt})
%
%
that the SoC-type architecture now outperforms the SBC-type architecture. Thus, the design of an embedded PS controller not only depends on its ability to provide solution trajectories within the time budget, but also on other factors such as available power and the system footprint.

Developing a proper architecture for an embedded PS controller is a nontrivial task. Significant care and caution are necessary to match the various aspects of the spectral algorithm to the computational units described in Fig.~\ref{fig:EMBEDDED_arch} while meeting practical requirements such as power, weight, volume etc. This is why the configuration of a ground system that may be needed for mission planning purposes is not necessarily the same as a flight system. These concepts are now embodied in commercially available embedded PS optimal control solutions, such as the KR8100 series available from Elissar Global, LLC~\cite{eg:press-release}.
Such systems are intended to help bridge the gap between PS optimal control theory and practice by offering an approach for rapid deployment of optimal control concepts to the field. A proliferation of these embedded controllers has the potential to revolutionize the application of optimal control concepts for emerging problems in aerospace and autonomous systems.

\section{Future Challenges}
\begin{figure}[h!]
      \centering\scalebox{1}
      {\includegraphics[scale=0.4]{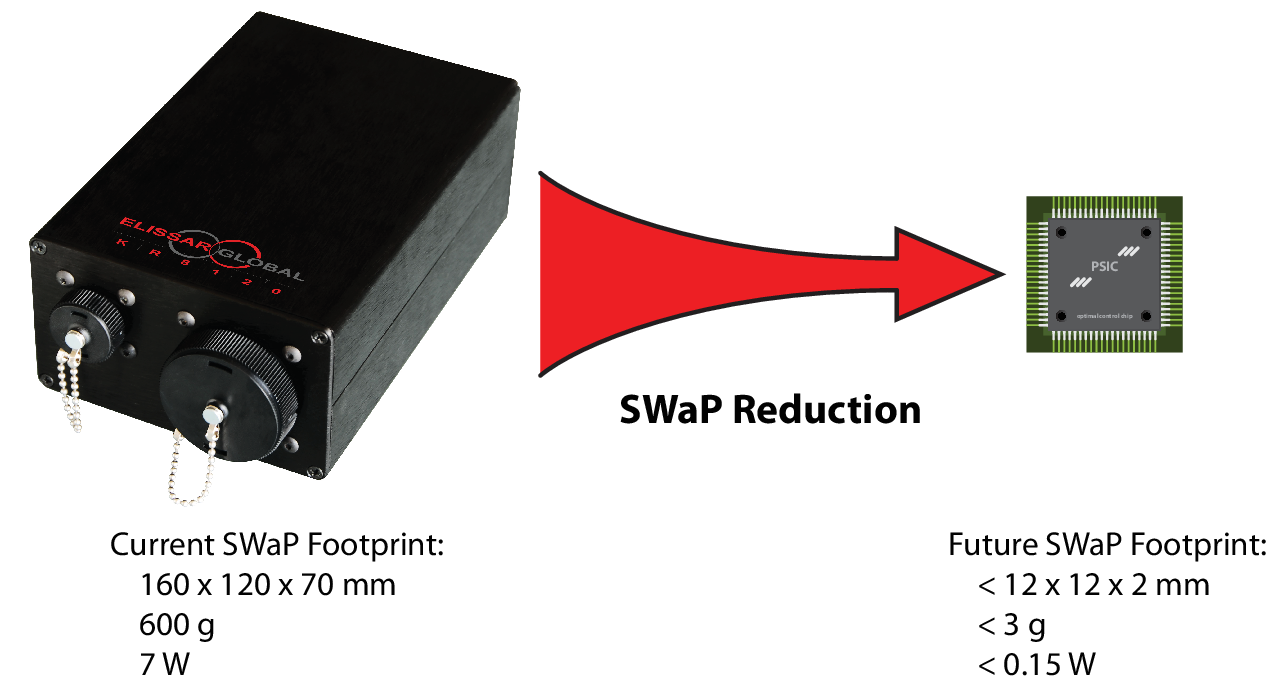}}
      \caption{\textsf{Future size, weight and power (SWaP) reduction challenges for embedded systems imposed by the space industry.}}
      \label{fig:KR8120}
   \end{figure}

PS optimal control started out as the Legendre pseudospectral method because of the elegant mathematical properties of Legendre polynomials.  The quest for improvement~\cite{arbGrid,Paul_Jacobi} over Legendre polynomials has not gone beyond the Chebyshev PS method; thus, these ``big two'' methods constitute the state of the art.  In this review paper, we have shown via unifying principles~\cite{RossGongBook,Fahroo_2008a,arbGrid} why all roads indeed lead to Legendre PS optimal control with the Chebyshev approach being its nearest neighbor. These unifying principles have generated a wealth of new ideas that have opened the door to a major new field of research in optimal control where theory, practice and software are intimately connected to flight success via the mathematical guarantees afforded by convergence theorems.  These theorems have been translated to algorithms implemented in DIDO and OTIS software packages, which account for the most widespread usage of PS optimal control techniques.
With success comes a growing demand for more.  The emerging challenges for PS control are vast and varied.  They range from a need to develop stronger theorems to all the way up to embedding PS control techniques in very small form factor PS integrated circuits; see Fig.~\ref{fig:KR8120}.
%
%
The need for the latter is driven by the space industry which requires form factors that are significantly smaller than present-day embedded systems~\cite{eg:press-release}.
This challenge is vast and requires a team of software engineers, numerical analysts, embedded system experts, optimization professionals and control practitioners to deliver a PS integrated circuit that can be plugged into a heritage system as a strap-on device as illustrated in Fig.~\ref{fig:PS-implement-block-diag}.

At the theoretical level, stronger convergence theorems are required; i.e. convergence theorems under weaker hypotheses. This need is driven, at least, in part because the computational results hold even when some of the assumptions are violated.  This suggests that the theorems are likely valid under a weaker set of hypotheses.

There are many applications (e.g mission planning~\cite{hybrid:jgcd}) wherein there is a need to generate PS controls faster than real time.  This need is driven by a demand to perform mission plans at a faster and faster pace. That is, in order to plan missions in near real time, each solution that is an element of the ``mission space'' needs to be computed extremely rapidly so that all solutions can be assembled to perform mission analysis. A key technical challenge here is to expand the scope of PS theory from standard optimal control theory to hybrid optimal control theory.  While major steps along  this direction have already been performed such as the development of a Hybrid Covector Mapping Principle~\cite{acc:hybrid}, the key challenge that remains in this area is in solving the accompanying graph problem. The graph problem is complicated by the fact that it is fundamentally coupled with the collection of standard optimal control problems that constitute the vertices of a digraph~\cite{hybrid:jgcd}. In this context, it can be convincingly argued that a nonlinear hybrid PS theory is one of the next grand challenges in control theory.



\section*{Acknowledgements}
We are deeply indebted to Profs. Fariba Fahroo, Qi Gong and Wei Kang, our long-time collaborators on PS optimal control theory and computation. It was Prof. Fahroo's initial work on costate estimation that led to the discovery of the covector mapping principle. Her ongoing efforts to further the ideas continues to inspire us all. Prof. Gong was instrumental in the proofs of the theorems on existence and consistency. We greatly appreciate his simple but powerful counter examples used in this paper to illustrate the subtle mathematical concepts. Prof. Kang led the way in proving the difficult convergence theorems and in estimating the rate of convergence. We are indebted to him for setting PS optimal control theory on a firm mathematical foundation. Last but not least, our heartfelt thanks goes to Dr. Nazareth Bedrossian for advancing an array of impressive implementations of PS control techniques onboard the International Space Station.








\end{document}